\documentclass[USenglish,onecolumn]{article}

\usepackage[utf8]{inputenc}%(only for the pdftex engine)
\usepackage[big]{dgruyter}
\usepackage[numbers]{natbib}
\usepackage{amsmath,amsfonts,mathtools}
\usepackage{multirow}
\usepackage{caption}
\usepackage{subfigure}
\usepackage{float}
\usepackage{booktabs, array, threeparttable}
\usepackage{xr}
\usepackage{xcolor}
\usepackage{hyperref}
\mathtoolsset{showonlyrefs=true}
\usepackage{enumitem}

\theoremstyle{plain}

\newtheorem{proposition}{Proposition}

\newtheorem{definition}{Definition}
\newtheorem{remark}{Remark}
\newtheorem{condition}{Condition}

\usepackage{etoolbox}
\usepackage{tabularx, booktabs, graphicx}

\makeatletter
\newcommand{\ostar}{\mathbin{\mathpalette\make@circled\star}}
\newcommand{\make@circled}[2]{%
  \ooalign{$\m@th#1\smallbigcirc{#1}$\cr\hidewidth$\m@th#1#2$\hidewidth\cr}%
}
\newcommand{\smallbigcirc}[1]{%
  \vcenter{\hbox{\scalebox{0.77778}{$\m@th#1\bigcirc$}}}%
}
\makeatother

\newcommand{\lt}{\left}
\newcommand{\rt}{\right}

\newcommand{\commenting}[1]{}

\newcommand{\Var}[1]{{\operatorname{Var}\left\{#1\right\}}}

\newcommand{\Cov}[2]{{\operatorname{Cov}\left\{#1,#2\right\}}}

\newcommand{\E}[1]{{\bbE\left\{#1\right\}}}
\newcommand{\Prob}[1]{{\bbP\left\{#1\right\}}}

\newcommand{\ind}[1]{\boldsymbol{1}\left\{#1\right\}}

\newcommand{\cnp}[2]{\begin{pmatrix}#1\\#2\end{pmatrix}} % choose p from n

\newcommand{\indep}{\perp \!\!\! \perp}

\ifdef{\see}{\renewcommand{\see}[1]{\text{ (#1)}}}{\newcommand{\see}[1]{\text{ (#1)}}}

\newcommand{\diag}[1]{\mathrm{Diag}\left\{#1\right\}}

\def\boxit#1{\vbox{\hrule\hbox{\vrule\kern6pt\vbox{\kern6pt#1\kern6pt}\kern6pt\vrule}\hrule}}

\newcolumntype{P}[1]{>{\centering\arraybackslash}p{#1}}
\newcolumntype{M}[1]{>{\centering\arraybackslash}m{#1}}
\newcolumntype{L}[1]{>{\raggedright\arraybackslash}m{#1}}

\DeclareMathOperator*{\argmin}{arg\,min}

\newcommand{\cA}{{\mathcal{A}}}

\newcommand{\cH}{{\mathcal{H}}}

\newcommand{\cL}{{\mathcal{L}}}

\newcommand{\cN}{{\mathcal{N}}}

\newcommand{\bbP}{{\mathbb{P}}}

\newcommand{\bbR}{{\mathbb{R}}}

\newcommand{\bbE}{{\mathbb{E}}}

\newcommand{\oa}{{\overline{a}}}
\newcommand{\ob}{{\overline{b}}}

\newcommand{\oY}{{\overline{Y}}}
\newcommand{\oX}{{\overline{X}}}

\newcommand{\halpha}{{\hat{\alpha}}}
\newcommand{\hgamma}{{\hat{\gamma}}}
\newcommand{\hbeta}{{\hat{\beta}}}
\newcommand{\heta}{{\hat{\eta}}}

\newcommand{\htau}{{\hat{\tau}}}

\newcommand{\hS}{{\hat{S}}}
\newcommand{\hV}{{\hat{V}}}
\newcommand{\hX}{{\hat{X}}}
\newcommand{\hY}{{\hat{Y}}}

\newcommand{\hDelta}{{\hat{\Delta}}}

\newcommand{\ttau}{{\tilde{\tau}}}

\newcommand{\scf}{{\textsc{f}}}
\newcommand{\scn}{{\textsc{n}}}
\newcommand{\scl}{{\textsc{l}}}

\newcommand{\scS}{{\textsc{S}}}
\newcommand{\scM}{{\textsc{M}}}

\def\Cov{\text{Cov}}

\usepackage[textsize=tiny, textwidth=4.5em]{todonotes}
\def\lxr{}

\begin{document}

  \articletype{Review Article {\hfill}Open Access}

  \author*[1]{Lei Shi}

\author[2]{Xinran Li}

  \affil[1]{University of California, Berkeley; E-mail: leishi@berkeley.edu}

  \affil[2]{The University of Chicago; E-mail: xinranli@uchicago.edu}

  \title{\huge Some theoretical foundations for the design and analysis of 
  randomized experiments
  }

  \runningtitle{A technical delving into
     randomized experiments} 

  %\subtitle{...}

  \begin{abstract}
  {
\citet{neyman1923application}'s seminal work 
% Neyman's seminal paper \cite{neyman1923application} 
in 1923 has been a milestone in statistics over the century, which has motivated many fundamental statistical concepts and methodology. In this review, 
we delve into \citet{neyman1923application}'s groundbreaking contribution 
and offer technical insights into the design and analysis of randomized experiments. We shall review the basic setup of completely randomized experiments and the classical approaches for inferring the average treatment effects.
We shall in particular review more efficient design and analysis of randomized experiments by utilizing pretreatment covariates, which move beyond Neyman's original work without involving any covariate. 
We then summarize several technical ingredients regarding randomizations and permutations that have been developed over the century, such as permutational central limit theorems and Berry--Esseen bounds,  
% (central limit theorems or CLTs and Berry--Esseen bounds or BEBs)  
and elaborate on how these technical results facilitate the understanding of randomized experiments. The discussion is also extended to other randomized experiments including rerandomization, stratified randomized experiments, matched pair experiments, cluster randomized experiments, etc. 
}
\end{abstract}
  \keywords{Causal inference, Permutation, Central limit theorem, Berry-Esseen bound, Potential outcome}
   \classification[MSC]{62K15, 62J05, 62G05}
 % \communicated{...}
 % \dedication{...}

  \journalname{Journal of Causal Inference}
\DOI{DOI}
  \startpage{1}
  \received{..}
  \revised{..}
  \accepted{..}

  \journalyear{2019}
  \journalvolume{1}
%  \journalissue{1}

\maketitle

\section{A review of the proposal in \citet{neyman1923application} and its influence}

Neyman's seminal paper \citep{neyman1923application} has been a cornerstone in the field of statistics over the last century. It has laid foundational principles that have significantly shaped multiple research areas such as causal inference, experimental design, and survey sampling. Its influence has been profound across a diverse range of applications, encompassing sectors like agriculture, economics, biomedical research, social science, and beyond.

The main purpose of  \citet{neyman1923application} is the analysis of field experiments conducted in order to compare a number of crop varieties. Suppose there are $m$ plots 
% (experimental units) 
and $\nu$ varieties.  
% (treatment arms). 
\citet{neyman1923application} introduced the notion of \textit{potential yield} of the $k$-th variety being applied to the $i$-th plot, which is denoted as $U_{ik}$, for $1\le i \le m$ and $1\le k \le \nu$; we use slightly different indices from Neyman to make them more intuitive.  In Neyman's framework, the quantities $\{U_{ik}\}$ are fixed but may be unknown. The number 
\begin{align*}
    a_k = \frac{1}{m} \sum_{i=1}^m U_{ik}
\end{align*}
is called ``the best estimate'' of the yield from the  $k$-th variety on the field, which is in fact an estimand representing the average yield in modern terminology. \citet{neyman1923application} then used an urn model as a thought experiment to depict the framework of sampling from a finite population. The $\nu$ types of varieties are treated as $\nu$ urns. Each urn contains $m$ balls, and each ball is associated with two labels: a plot label indexing the plots and a yield label indicating the unknown potential yields on the plots for each of the varieties. Specifically, in the $k$-th urn, there are $m$ balls with yield labels
\begin{align}
    U_{1k}, \dots, U_{ik}, \dots, U_{mk}. 
\end{align}
Also, the urns have the property that ``if one ball is taken from one of them, the balls having the same plot label disappear from all the other urns''. Then from each urn, a number of balls are drawn without replacement. With this model, Neyman studied the properties (in particular the means and variances) of the sample 
averages across 
% treated and control arms 
{all varieties}
as well as their difference under the 
% permutational 
{randomization}
distribution. This marks the pioneering effort for studying the difference-in-means estimator in modern terminology. Notably, he was able to ``determine empirically that the difference of partial averages of the plots sampled shows a fair agreement with the Gaussian law distribution''. This corrects the ``common misunderstanding'' at that time that inference can be performed only if the yields from different plots follow the Gaussian law. Combined with a {conservative} variance estimation strategy, he suggested 
a confidence interval for 
% inference 
{the true difference between two varieties}
based on normal approximation. % Common misunderstanding consists in the unjustified assertion that probability theory can be applied to solve problems similar to the one discussed only if the yields from the different plots follow the Gaussian law 

\citet{neyman1923application} offered a series of groundbreaking and foundational insights. Below we outline three key facets of \citet{neyman1923application}'s contributions.

The first contribution is the introduction of the potential outcome model. This model has since become a standard framework for illustrating possible experimental outcomes, as referenced in works such as \cite{pitman1937significance, welch1937z, kempthorne52, kempthorne1955randomization}. The potential outcome paradigm serves as an impeccable model for a discussion in causation within randomized experiments. Within this framework, researchers pose and address causal questions by analyzing causal effects
% , which 
that are defined as comparisons between potential outcomes, 
which represent various hypothetical scenarios or states of the world. This framework also elegantly facilitates the representation of interference between units \citep{hudgens2008toward, tchetgen2012causal, zigler2021bipartite}, the prolonged impacts of interventions \citep{liu2022practical, sjolander2016carryover, imai2023matching}, and the causal analyses involving post-treatment variables such as instrumental variables \citep{angrist1996identification} and mediation \citep{vanderweele2015explanation, vanderweele2016mediation}. Moreover, the importance of potential outcomes transcends experimental settings and is also profound in observational studies, as highlighted by \citet{rubin1974estimating}.

The second contribution of \citet{neyman1923application} lies in that it further highlights the importance of physical randomization or random selection when conducting experiments or performing sampling. Randomization has been in the air since 1920s, as commented by \citet{rubin1990comment} citing \citet{student1923testing} and \citet{fisher1923studies} as references. \citet{neyman1923application} contributed to the randomization world by 
% appling 
introducing the potential outcome model and describing a finite population inference framework for randomization. Within this framework, potential outcomes are viewed as fixed, and physical randomization emerges as the ``reasoned basis'' \citep{fisher1935design} for facilitating statistical testing and estimation \citep{kempthorne52, rosenbaum2002observational, hinkelmann2007design, imbens15}. Moreover, the proposal of sampling without replacement also inspires the pursuit of the parallels and linkages between survey sampling and randomized experiments \citep{splawa1925contributions, neyman1934two, neyman1935statistical, fienberg1996reconsidering}.

The third contribution of \citet{neyman1923application} centers on the repeated sampling properties of statistics over their non-null randomization distribution. This viewpoint offers a new perspective on randomization-based or design-based inference, distinguishing it from Fisher's focus on the sharp null hypothesis of no causal effects
for any units and finite-sample exact p-values \citep{rubin1990comment}.   \citet{neyman1923application}  recognized from an empirical perspective that the asymptotic normality holds under the 
% permutational 
described sampling scheme, without requiring the outcomes to come from a Gaussian law. Moreover, he proposed 
to estimate the variance of an estimator conservatively in expectation, which can further lead
% a variance estimation strategy that led 
to a conservative confidence interval. These efforts built up the foundation for 
% asymptotically valid 
large-sample randomization-based inference in finite populations. 

Building upon the pioneering contribution of \citet{neyman1923application} in randomization-based inference, there have been many new developments in the design and analysis of 
% {\lxr randomized experiments}. 
randomized experiments.
In the following sections, we shall first review the basic setup of completely randomized experiments (CREs) and the classical approaches for analysis. We then present several technical ingredients regarding randomizations and permutations, such as central limit theorems (CLTs) and Berry--Esseen bounds (BEBs),  
% (central limit theorems or CLTs and Berry--Esseen bounds or BEBs) 
that were developed over the century,  
and elaborate on how these results enhance and expand our understanding of the design and analysis of CREs. We also extend the discussion to other randomized experiments and permutation-related technical tools.

\textit{Notations.} We summarise a set of notations for the whole paper. For an integer $N$, we use $[N]$ to denote the set of integers $\{1,\dots, N\}$. For two positive semidefinite matrices $V_1$ and $V_2$, we use $V_1 \gtrsim V_2$ to indicate that $V_1$ dominates $V_2$, in the sense that $V_1 - V_2$ is positive semidefinite. For a random sequence $\{X_N\}_{N=1}^\infty$, we write $X_N \rightsquigarrow \cL$ if $X_N$ converges weakly to the distribution $\cL$ as $N\to\infty$, and $X_N \rightsquigarrow L_N$ if $X_N$ and $L_N$ converge weakly to the same distribution. When $X_N$ converges in probability, we use $\mathrm{plim}_{N\to\infty}X_N$ to denote its probabilistic limit.

\section{Design and analysis of completely randomized experiments}\label{sec:design-analysis}

In this section, we introduce the basic setup for the design and analysis of CREs. Section \ref{sec:basic-cre} discusses the setup of a simple treatment-control CRE as well as strategies for estimation and inference. The results are extended to a  more general multi-level CRE. We then consider more efficient design and analysis of randomized experiments by incorporating pretreatment covariates. 
In particular, Section \ref{sec:cov-adj} presents several covariate-adjusted estimators, and Section \ref{sec:rerandomization} discusses rerandomization. 

% Section \ref{sec:cov-adj} presents several covariate-adjusted estimators. Section \ref{sec:rerandomization} discusses rerandomization, which essentially conducts covariate adjustment in the design stage. 

\subsection{Basic design and analysis of completely randomized experiments}\label{sec:basic-cre}

\subsubsection{Treatment-control completely randomized experiment}

We start 
by considering 
a treatment-control CRE 
that enrolls $N$ units, 
with $N_1$ units in the treatment arm and $N_0$ in the control arm. Let $Z_i$ denote the treatment assignment indicator for the $i$-th unit, for $1\le i \le N$. The treatment assignment status for the entire experiment is vectorized as   $Z = (Z_1,\dots, Z_N)$. Under complete randomization,
\begin{align}
    \Prob{Z = (z_1,\dots, z_N) } = 1/\cnp{N}{N_1}, ~\text{ for any } (z_1,\dots, z_N) \in \{0,1\}^N \text{ with }\sum_{i=1}^N z_i = N_1, ~\sum_{i=1}^N (1-z_i) = N_0.
\end{align}
The potential outcomes for the $i$-th unit are $(Y_i(1), Y_i(0))$. 
This is essentially a special case of \citet{neyman1923application}'s setup with two interventional arms. 
The more general notions of experimental units, treatment/control arms, and potential outcomes presented here correspond to \citet{neyman1923application}'s notions of plots, varieties, and potential yields.
% which corresponds to Neyman's notion of potential yields under the hypothetical interventions 1 and 0. 

\citet{rubin2005causal} called the 
% following 
$N\times 2$ matrix 
of potential outcomes in Table \ref{tab:science-tab}
as the Science Table. 
\begin{table}[!ht]
    \centering
    \caption{Science Table for Treatment-control CRE}
    \label{tab:science-tab}
    \begin{tabular}{ccc}
        $i$ & $Y_i(0)$ & $Y_i(1)$ \\
    \midrule
        $1$ & $Y_1(0)$ & $Y_1(1)$ \\
        $\vdots$ & $\vdots$ & $\vdots$ \\
        $N$ & $Y_N(0)$ & $Y_N(1)$ \\
    \end{tabular}
\end{table}
The observed outcome for the $i$-th unit is $Y_i = Z_iY_i(1) + (1-Z_i)Y_i(0)$, one of the two potential outcomes. Importantly, the potential outcomes are fixed and the randomness comes merely from the random allocation of the treatment, reflected by the random vector $Z$. \citet[][Chapter 9]{scheffe1959analysis} called it the {\it randomization model}. Under this model, it is conventional to call the resulting inference as {\it randomization-based inference}, {\it design-based inference} or {\it finite population inference}.  It has become increasingly popular in both theory and practice \citep[e.g.][]{kempthorne52, copas1973randomization, robins1988confidence, rosenbaum2002observational, hinkelmann07, freedman2008Aregression, freedman2008Bregression, lin2013agnostic, dasgupta2015causal, imbens15, ATHEY201773, fogarty2018regression, guo2021generalized}.
Define further the following finite-population mean and variance of potential outcomes for each arm, which are essentially summaries of the science table in Table \ref{tab:science-tab}:
\begin{gather}
    \oY(0) = \frac{1}{N}\sum_{i=1}^N Y_i(0), \quad
    \oY(1) = \frac{1}{N}\sum_{i=1}^N Y_i(1);\label{eqn:means}\\
    S^2(0) = \frac{1}{N-1}\sum_{i=1}^N (Y_i(0) - \oY(0))^2, \quad S^2(1) = \frac{1}{N-1}\sum_{i=1}^N (Y_i(1) - \oY(1))^2. \label{eqn:vars}
\end{gather}

Under the potential outcome framework,
the $i$-th unit has individual treatment effect $\tau_i = Y_i(1) - Y_i(0)$, for $1\le i \le N$. The average treatment effect (ATE) over all units is then defined as:
\begin{align}
    \tau = \frac{1}{N}\sum_{i=1}^N \tau_i = \oY(1) - \oY(0).
\end{align}
% \citet{neyman1923application} studied the problem of testing the following null hypothesis:
% \begin{align}\label{eqn:weak-null}
%     H_{0\textsc{N}}: \tau = \overline{Y}(1) - 
%     \overline{Y}(0) = 0,
% \end{align}
% which is also called the weak null hypothesis \citep{wu2021randomization}.
% For \eqref{eqn:weak-null}, a straightforward testing statistic is 
\citet{neyman1923application} proposed to estimate the ATE $\tau$ by 
the difference-in-means estimator: 
\begin{align}\label{eqn:htau}
    \htau = \hY(1) - \hY(0), \text{ where } \hY(z) =  \frac{1}{N_z}\sum_{i=1}^N Y_i\ind{Z_i = z}, \text{ for } z=0,1.
\end{align}
% Under \eqref{eqn:weak-null}, 
% \citet{neyman1923application} 
He proved that $\htau$ is an unbiased estimator for $\tau$, i.e., $\E{\htau} = \oY(1) - \oY(0) = \tau$,  
with true variance
\begin{align*}
    \Var{\htau} = \frac{1}{N_1} S^2(1) + \frac{1}{N_0} S^2(0) - \frac{1}{N} S^2(\tau), 
\end{align*}
where the variances $S^2(0)$ and $S^2(1)$ are defined in \eqref{eqn:vars} and $S^2(\tau)$ is the variance of the individual treatment effects
\begin{align}\label{eqn:Stau}
    S^2(\tau) = \frac{1}{N-1}\sum_{i=1}^N (\tau_i - \tau)^2.
\end{align}

Due to the fact that we are never able to jointly observe the two potential outcomes for any unit, the variance of individual effects in \eqref{eqn:Stau} is generally not estimable based on the observed data. \citet{neyman1923application} proposed the following variance estimator:
\begin{align}\label{eqn:vN}
    \hV = \frac{1}{N_1} \hS^2(1) + \frac{1}{N_0} \hS^2(0), \text{ where } \hS^2(z) = \frac{1}{N_z-1}\sum_{i=1}^N (Y_i - \hY(z))^2\ind{Z_i = z}.
\end{align}
which essentially circumvents the problem by dropping the unestimable component regarding $S^2(\tau)$. 
The variance estimator in 
\eqref{eqn:vN} has expectation
\begin{align}\label{eqn:EhvN}
    \E{\hV} = \frac{1}{N_1}S^2(1) + \frac{1}{N_0}S^2(0) \ge \Var{\htau},
\end{align}
which suggests that $\hV$ is in general not unbiased but conservative. A level-$\alpha$ confidence interval is then given by
\begin{align}\label{eqn:CI}
\lt[\htau - z_{{\alpha}/{2}}\sqrt{\hV}, ~ \htau + z_{{\alpha}/{2}}\sqrt{\hV}\rt],
\end{align}
where $z_{\alpha/2}$ is the $\alpha/2$ upper quantile of a standard normal distribution. In Sections \ref{sec:PCLT} and \ref{sec:PBEB} we will discuss more technical results for the asymptotic validity of the confidence interval in \eqref{eqn:CI}.

\begin{remark}\label{rmk:remark_neyman_fisher}
\citet{neyman1923application}'s approach can also be used to test the following null hypothesis:
\begin{align}\label{eqn:weak-null}
    H_{0\textup{N}}: \tau = \overline{Y}(1) - 
    \overline{Y}(0) = 0,
\end{align}
which is often called the weak null hypothesis \citep{wu2021randomization}.
In contrast, 
\citet{fisher1935design} proposed to test the following null hypothesis:
\begin{align}\label{eqn:strong-null}
    H_{0\textup{F}}: Y_i(1) = Y_i(0) \text{ for all units } i = 1,\dots,N, 
\end{align}
which is called the sharp null hypothesis by \citet{rubin1980randomization} or the strong null hypothesis by \citet{wu2021randomization}.
The Fisherian perspective is fundamentally different, as it focuses on testing the hypothesis of no causal effects for any units whatsoever, whereas the Neymanian perspective focuses on testing no average causal effect \citep{ding2017paradox}. Obviously, Fisher's null implies Neyman's null, but either of them can be practically relevant depending on the application. 
Under \eqref{eqn:strong-null}, one can impute the unobserved potential outcomes and perform  Fisher's randomization test to deliver finite sample exact inference \citep{fisher1935design}. 
Fisher's test has the advantage of being finite-sample valid, while Neyman's requires large-sample approximation. 
% However, the tested null hypothesis by Fisher is fundamentally different from Neyman's.  Intuitively, Fisher is testing whether a treatment is irrelevant for the outcome, while Neyman focuses on the average effects of the treatment; 
% I think the key advantage of Fisher’s null is not just being finite-sample valid (in fact, from my experiences, today’s world is less and less dependent on using a small sample to make critical decisions, so one may argue finite-sample validity is good to have but not critical). Many people, including Fisher, would argue no average effect is fundamentally different from no effect. Sharp null is testing “irrelevance,” which is quite different from no average effect.
% more restrictive than that by Neyman, since the sharp null requires speculation of all individual treatment effects.
Nevertheless, Neyman's asymptotic results can also help ease the computation for Fisher's null.
We refer interested readers to refs \cite{ding2018randomization, zhao2021covariate, cohen2022gaussian} for a unification of both perspectives, 
and to refs \cite{Rosenbaum:2001, rigdon2015exact, LD2016binary, CDLM21quantile, SL2023, chen2023role} for extending FRT to nonsharp null hypotheses.

% % Fisher's sharp null hypothesis is closely related to \citet{neyman1923application}'s proposal, {\lxr although it is more restrictive in the sense that it requires speculation of all individual treatment effects}. 
% A unification of both perspectives is discussed in \citep{ding2018randomization, zhao2021covariate, cohen2022gaussian}.
% As a side note, Neyman's asymptotic results can also help ease the computation for Fisher's null. 

\end{remark}

\begin{remark}\label{rmk:regression_simple}
For analysis, practitioners usually prefer regression-based inference for the average causal effect. The standard approach is to run the ordinary least squares (OLS) of the
outcomes on the treatment indicators with an intercept:
\begin{align}\label{eq:ols}
    (\hgamma, \htau) = \argmin_{\gamma, \tau\in\bbR} \sum_{i=1}^N (Y_i - \gamma - Z_i \tau)^2. 
\end{align}
As implicitly written in \eqref{eq:ols}, the point estimator from the OLS for the treatment effect is identical to the difference-in-means estimator in \eqref{eqn:htau}. 
However, the usual variance estimation based on the OLS usually fails (in the sense of either underestimating or overestimating the truth by possibly a quite large factor), 
due to heteroskedasticity in potential outcomes \citep{freedman2008Aregression}. 
More concretely, the OLS-based variance estimator is 
\begin{align}
\hat{V}_{\mathrm{OLS}} & =\frac{N\left(N_1-1\right)}{(N-2) N_1 N_0} \hat{S}^2(1)+\frac{N\left(N_0-1\right)}{(N-2) N_1 N_0} \hat{S}^2(0) 
\approx \frac{\hat{S}^2(1)}{N_0}+\frac{\hat{S}^2(0)}{N_1},
\end{align}
which can be very different from \eqref{eqn:hVN} if the numbers of units or the sample variances of observed outcomes differ a lot between the two arms. 
Instead, one can use the Eicker-Huber-White (EHW) variance estimator to obtain a robust estimation:
\begin{align}\label{eqn:EHW}
    \hV_{\textup{EHW}} = \frac{\hS^2(1)}{N_1}\frac{N_1 - 1}{N_1} + \frac{\hS^2(0)}{N_0}\frac{N_0 - 1}{N_0},
\end{align}
which is asymptotically equivalent to $\hV$ in \eqref{eqn:vN}. Alternatively, the so-called HC2 variant of the EHW robust variance estimator is identical to $ \hV $.
See Chapter 4 of  ref. \cite{ding2023first} for more detailed discussion on regression-based analyses for the average treatment effect. 
\end{remark}

% \textbf{Multi-level CRE.} 
\subsubsection{Multi-level completely randomized experiments}\label{sec:multi_CRE}

Many efforts have been devoted to extending the treatment-control CRE to multi-level scenarios, which caters for many 
% realistic questions 
practical problems 
and designs such as (fractional) factorial experiments \citep{wu2011experiments, dasgupta2015causal}, conjoint analysis \citep{hainmueller2014causal, hainmueller2015hidden}, partially nested experiment \citep{bauer2008evaluating, hallfors2006efficacy}, 
sampling-based randomized experiments \citep{BD20, YQL2021}, etc.

In a multi-level randomized experiment, there are $N$ units and $Q$ treatment arms, where the number of units under treatment $q$ equals $N_q$, with $\sum_{q=1}^Q N_q= N$. Corresponding to treatment level $q $, unit $i$ has the potential outcome $Y_i(q)$, where $i=1, \ldots, N$ and $q = 1, \ldots, Q$; see the Science Table in Table \ref{tab:science-tab-multi-level}. Despite its simplicity, the multi-level CRE has been widely used in practice and has generated rich theoretical results. Definition \ref{def:cr} below characterizes the joint distribution of $Z =  (Z_1, \ldots, Z_N)$ under complete randomization, where $Z_i \in \{1, \ldots, Q\}$ is the treatment indicator for unit $i$. 

\begin{definition}[Complete randomization]\label{def:cr}
Fix integers $N_1,\ldots, N_Q$ with $\sum_{q=1}^Q N_q= N$. The treatment vector $Z$ is 
uniformly distributed over $\mathcal{Z} \equiv \{z\in \{1,2,\ldots,Q\}^N: \sum_{i=1}^N \ind{z_i = q} = N_q, \text{ for } 1\le q\le Q \}$. 
% uniform over all possible values. 
\end{definition}

Mathematically, Definition \ref{def:cr} implies that $ \bbP( {Z} = {z} ) = N_1 ! \cdots N_Q! / N! $ for all possible values of $z$ in $\mathcal{Z}$. 
% ${z} = (z_1, \ldots, z_N ) $. 
Computationally,  Definition \ref{def:cr} implies that ${Z}$ is from a random permutation of $N_1$ $1$'s, $N_2$ $2$'s, $\ldots $, $N_Q$ $Q$'s. The observed outcome is $Y_i = \sum_{q=1}^Q Y_i(q)\ind{Z_i = q}$ for each unit $i$. 

Similar to the two-arm setting discussed in Section \ref{sec:basic-cre}, in \citet{neyman1923application}'s framework, all potential outcomes are fixed and only the treatment indicators are random according to Definition \ref{def:cr}. 

\begin{table}[!ht]
    \centering
    \caption{Science Table for Multi-level CRE}
    \label{tab:science-tab-multi-level}
    \begin{tabular}{ccccc}
        $i$ & $Y_i(1)$ & $Y_i(2)$ & $\dots$ & $Y_i(Q)$ \\
    \midrule
        $1$ & $Y_1(1)$ & $Y_1(2)$ & $\dots$ & $Y_1(Q)$ \\
        $\vdots$ & $\vdots$ & $\vdots$ & $\ddots$  & $\vdots$\\
        $N$ & $Y_N(1)$ & $Y_N(2)$ & $\dots$ & $Y_N(Q)$ \\
    \end{tabular}
\end{table}

We consider a general contrast matrix $F\in\bbR^{Q\times H}$ of full column rank, i.e, $F^\top 1_Q = 0_H$ and $\operatorname{rank}(F) = H$, 
and a set of individual treatment effects defined as the linear contrasts of the potential outcomes: 
% Neyman's weak null hypothesis is also extended to the multi-level experiments  \citep{wu2021randomization}. 
% Let $F\in\bbR^{Q\times H}$ be a full column rank contrast matrix, i.e. $F^\top 1_Q = 0_H$.  
% Then a general set of individual treatment effects can be defined as 
\begin{align}\label{eq:tau_i_F}
    \tau_i = F^\top Y_i(\cdot),
\end{align}
where $Y_i(\cdot)$ is the vectorized potential outcomes for unit $i$:
\begin{align}
    Y_i(\cdot) = (Y_i(1), \dots, Y_i(Q))^\top. 
\end{align}
The average effect is defined as 
\begin{align}\label{eqn:ATE-Q}
    \tau = \frac{1}{N}\sum_{i=1}^N \tau_i = F^\top \oY(\cdot),
\end{align}
where $\oY(\cdot)$ is the vectorized average potential outcomes:
\begin{align*}
    \oY(\cdot) = \frac{1}{N}\sum_{i=1}^N Y_i(\cdot) = (\oY(1), \cdots, \oY(Q))^\top. 
\end{align*}
% Neyman's weak null can be defined as follows:
% \begin{align*}
%     H_{0\textsc{N}}: \tau = F^\top \oY(\cdot) = 0.
% \end{align*}
When $Q=2$ and $F = (1, -1)^\top$, 
$\tau$ in \eqref{eqn:ATE-Q} reduces to
% corresponds to 
the ATE in the treatment-control setting. Moreover, 
we can estimate $\tau$
% using the plug-in estimator, which can be viewed as a 
by the following 
generalization of the difference-in-means estimator:
\begin{align}\label{eqn:htau-Q}
    \htau = F^\top \hY(\cdot),
\end{align}
where $\hY(\cdot) = (\hY(1), \dots, \hY(Q))^\top$ is the vectorized sample averages of observed outcomes for all treatment arms, with $\hY(q) = {N_q^{-1}} \sum_{i=1}^N Y_i \ind{Z_i = q}$. 
The estimator in 
\eqref{eqn:htau-Q} has variance \citep{li2017general}
\begin{align}\label{eqn:var-htau-Q}
    \Var{\htau} = F^\top \diag{\frac{1}{N_q} S(q,q)}_{q=1}^Q F - \frac{1}{N} F^\top S F,
\end{align}
where $S \in \bbR^{Q\times Q}$ is a covariance matrix for the potential outcomes with the $(q,q')$-th entry given by 
\begin{align}\label{eq:S_qq}
    S(q,q') = \frac{1}{N - 1}\sum_{i=1}^N (Y_i(q) - \oY(q))(Y_i(q') - \oY(q')), \quad q,q'=1,\dots,Q,
\end{align}
and $F^\top S F$ is essentially the finite population covariance of the individual effects $\tau_i$'s in \eqref{eq:tau_i_F}. 
A variance estimator for \eqref{eqn:var-htau-Q} is 
\begin{align}\label{eqn:hVN}
\hV = F^\top \diag{\frac{1}{N_q} \hS(q,q)}_{q=1}^Q F, 
\end{align}
where $\hS(q,q)$ is the sample variance within treatment level $q$:
\begin{align}\label{eq:S_hat_qq}
    \hS(q,q) = \frac{1}{N_q - 1}\sum_{i=1}^N (Y_i - \hY(q))^2 \ind{Z_i = q}. 
\end{align}
Using \eqref{eqn:htau-Q} and \eqref{eqn:hVN}, a Wald-type confidence region for $\tau$ is given by 
\begin{align}\label{eqn:CR}
    \lt\{\tau: (\htau - \tau)^\top \hV^{-1} (\htau - \tau) \le q_{H, \alpha}\rt\},
\end{align}
where $q_{H, \alpha}$ is the upper-$\alpha$ quantile of the $\chi^2_H$ distribution.  \eqref{eqn:CR} can be proved to be asymptotically valid under mild regularity conditions. More details are deferred to Sections \ref{sec:PCLT} and \ref{sec:PBEB}. 

Below we give two remarks in parallel with Remarks \ref{rmk:remark_neyman_fisher} and \ref{rmk:regression_simple}. 
First, 
% As a final remark,  
Fisher's randomization test (FRT) has also been a popular tool for analyzing multiple-level randomized experiments, which can be used to test sharp nulls and deliver finite-sample exact p-values \citep{imbens15}.
% There have been many efforts to use FRTs to test weak null hypotheses. 
See also \citep{wu2021randomization, zhao2021covariate, cohen2022gaussian} for the unification of Neyman's and Fisher's approaches. 
Second, similar to the treatment-control case, we can perform analysis with the regression-based approach. \citet{zhao2023covariate}
studied general regression-based analyses in multi-level experiments. % As a special case, \citet{lu2023design, zhao2022regression}, among others, studied regression methods for factorial experiments. 

% {\color{red} a regression-based analysis framework... \citet{zhao2022regression} studied the use of weighted least squares (WLS) for analyzing }

\subsection{Covariate adjustment}\label{sec:cov-adj}

In many randomized experiments, there are pre-treatment covariates $X_1,\dots, X_N$ for the $N$ units, where $X_i$'s are encoded as vectors in $\bbR^p$. Covariate adjustment has become a standard approach for analyzing randomized experiments and has been widely adopted in many fields. As one example, in 2023, U.S. Food and Drug Administration issued the final guidance on \textit{Adjusting for Covariates in Randomized Clinical Trials for Drugs and Biological Products Final Guidance for Industry}. This guidance describes the agency’s current recommendations regarding adjusting for covariates in the statistical analysis of randomized clinical trials in drug and biological product development programs.
A natural question is 
{\lxr how to optimally adjust covariates for inference?}
% whether adjusting for covariates is helpful for inference. 
The problem is nontrivial in several aspects: (i) the true relation between outcomes and covariates is usually unknown;  (ii) the potential outcomes under different treatment levels are in general heterogeneous. Many research works explored covariate adjustment from both practical and theoretical perspectives.  It has become a standard practice to use a model-assisted method for covariate adjustment to gain efficiency for inference 
while being robust to model misspecification \citep{lin2013agnostic}.
% when the underlying model is misspecified.

% \textbf{Fisher's analysis of covariance (ANCOVA)}.

\subsubsection{Fisher's analysis of covariance (ANCOVA)}
Historically, \citet{fisher1925statistical} proposed to use the analysis of covariance (ANCOVA) to improve estimation efficiency. This remains a standard strategy in
many fields. He suggested running the OLS of $Y_i$ on $(1, Z_i
, X_i)$ and using
the coefficient of $Z_i$ as an estimator for $\tau$. Mathematically, let 
$\oX$ be the mean of the covariates: $\oX = {N}^{-1}\sum_{i=1}^N X_i$. 
Fisher’s ANCOVA estimator $\htau$ is given by the following OLS output:
\begin{align}\label{eqn:fisher-ancova}
    (\htau_\scf,\halpha_\scf,\hbeta_\scf) = \argmin_{\alpha,\tau\in\bbR, \beta\in\bbR^p} 
    % \frac{1}{2N}
    \sum_{i=1}^N \{Y_i - \alpha - Z_i \tau - (X_i - \oX)^\top \beta\}^2,
\end{align}
noting that the centering of covariates in \eqref{eqn:fisher-ancova}
% is 
% not necessary and 
will not affect the OLS estimator $\htau_\scf$.

\citet{freedman2008Aregression, freedman2008Bregression} studied Fisher's ANCOVA estimator under the CRE. He showed that $\hat{\tau}_F$ can be biased in the finite sample, but is consistent for the true average effect as the sample size goes to infinity. Moreover, he showed some negative results for Fisher's ANCOVA estimator. First, the asymptotic variance of $\htau_\scf$ can be even larger than the simple difference-in-means estimator $\htau$ without adjusting any covariates. Second, the standard error estimator from OLS can underestimate the true standard error of $\htau_\scf$ under the CRE.

\subsubsection{Lin's estimator}\label{sec:Lin}
% \textbf{Lin's estimator.} 
In response to Freedman's negative findings, \citet{lin2013agnostic} proposed a remedy, which is called ``Lin's estimator'' nowadays. Concretely speaking, he proposed to run OLS of $Y_i$ on $Z_i$ and $X_i$ as well as their interaction term $Z_i \times X_i$:
\begin{align}\label{eqn:lin}
    (\htau_\scl, \halpha_\scl, \hbeta_\scl, \heta_\scl) = \argmin_{\alpha,\tau\in\bbR, \beta,\eta\in\bbR^p} \frac{1}{2N}\sum_{i=1}^N \{Y_i - \alpha - Z_i \tau - (X_i - \oX)^\top \beta - Z_i \times (X_i - \oX)^\top \eta\}^2.
\end{align}
Importantly, unlike \eqref{eqn:fisher-ancova}, the centering of covariates here is critical since it will affect the OLS estimator $\htau_\scl$. 

% {\color{red}
% [R1: In the paragraph below, the statement ”it is consistent when the sample size N goes to infinity” is
% made. However, since Fisher’s ANCOVA estimator is also consistent, consistency should not be
% considered a benefit. On the contrary, both estimators are biased in finite samples, which is a
% drawback.]
% }

Lin's estimator is also consistent when the sample size $N$ goes to infinity. Moreover, it enjoys several benefits. First, the asymptotic variance of $\htau_\scl$ is no larger than that of both $\htau$ and $\htau_\scf$. Second, the 
% Eicker-Huber-White (EHW) 
EHW variance estimator for 
% \eqref{eqn:fisher-ancova} and 
\eqref{eqn:lin} is asymptotically conservative for the true variance of 
% $\htau_\scf$ and 
$\htau_\scl$.
% respectively.
As a side note, the EHW standard error estimator for \eqref{eqn:fisher-ancova} is also asymptotically conservative for the true variance of $\htau_\scf$. See \citet{lin2013agnostic} for a more formal presentation of the theoretical results.

Besides the regression proposal, a second perspective for understanding Lin's estimator is based on minimizing the true or estimated variance of linearly adjusted estimators \citep{cochran1977sampling}. 
% Here we take the estimated variance perspective. 
Consider the following class of linearly covariate-adjusted estimators:
\begin{align}\label{eqn:linear-estimators}
    \htau(\beta_1, \beta_0) &= N_1^{-1} \sum_{i=1}^N Z_i\{Y_i - (X_i-\oX)^\top \beta_1\} - N_0^{-1} \sum_{i=1}^N (1-Z_i)\{Y_i - (X_i-\oX)^\top \beta_0\} \\
    &= \{\hY(1) - (\hX(1)-\oX)^\top\beta_1\} - \{\hY(0) - (\hX(0)-\oX)^\top\beta_0\} \\
    % = \hgamma_1 - \hgamma_0.
    & = \htau - \delta^\top \htau_X, 
\end{align} 
where $\hX(1)$ and $\hX(0)$ denote the averages of covariates in treatment and control groups,
$\htau_X \equiv \hX(1)-\hX(0)$ denotes the difference-in-means of covariates,
and $\delta = N_0/N \cdot \beta_1 + N_1/N \cdot \beta_0$ is a weighted average of the two linear adjustment coefficients. 
Obviously, the true variance of the covariate-adjusted estimator in \eqref{eqn:linear-estimators} is minimized when $\delta$ is the least squares coefficient from regressing the difference-in-means estimator $\hat{\tau}$ on the difference-in-means of covariates $\hat{\tau}_X$ under the CRE. 
\citet{li2017general} showed that this is further achieved when $\beta_1$ and $\beta_0$ are the least squares coefficients from projecting the treatment and control potential outcomes on covariates, respectively. 
Moreover, since the potential outcomes cannot be fully observed, we can estimate the least squares coefficients by their sample analogues $\hat{\beta}_1$ and $\hat{\beta}_0$, which are the least squares coefficients from the linear projection of observed outcomes on covariates in treatment and control groups, respectively. 
The resulting covariate-adjusted estimator 
$\htau - \hat{\delta}^\top \htau_X$ with $\hat{\delta} = N_0/N \cdot \hbeta_1 + N_1/N \cdot \hbeta_0$ is actually identical to Lin's estimator. 

We consider then the estimated variance for the covariate-adjusted estimator in \eqref{eqn:linear-estimators}.
We can essentially view the covariate-adjusted estimator as the difference-in-means estimator but with the adjusted potential outcomes $Y_i(1) - (X_i-\oX)^\top \beta_1$ and $Y_i(0) - (X_i-\oX)^\top \beta_0$. 
From the discussion in Section \ref{sec:basic-cre}, a conservative variance estimator for \eqref{eqn:linear-estimators} can be 
% is given by
\begin{align}\label{eqn:hvN-beta}
    \hV(\beta_1, \beta_0) 
    &= \{N_1(N_1-1)\}^{-1} \sum_{i=1}^N Z_i\{Y_i - \hgamma_1 - (X_i-\oX)^\top \beta_1\}^2 \\
    &+ 
    \{N_0(N_0-1)\}^{-1} \sum_{i=1}^N (1-Z_i)\{Y_i - \hgamma_0 - (X_i-\oX)^\top \beta_0\}^2,
\end{align}
where $\hgamma_1$ and $\hgamma_0$ are the sample mean of the adjusted outcomes for the treatment and control arm, respectively:
\begin{align}
    \hgamma_1 = 
    \frac{1}{N_1}\sum_{i=1}^N
    Z_i\{Y_i - (X_i - \oX)^\top \beta_1\}, 
    \quad 
    \hgamma_0 = \frac{1}{N_0}\sum_{i=1}^N
    (1-Z_i)\{Y_i - (X_i - \oX)^\top \beta_0\}.
\end{align}
This formulation suggests choosing $\beta_1$ and $\beta_0$ to minimize the variance estimator $\hV(\beta_1, \beta_0)$ to obtain a plug-in estimator for $\beta_1$ and $\beta_0$, which is equivalent to solving the following two regression problems for treated and control groups, respectively, with intercept terms $\gamma_1$ and $\gamma_0$ \citep{ding2023first}:
\begin{align}\label{eqn:two-regressions}
    \min_{\gamma_1, \beta_1} \sum_{i=1}^N Z_i\{Y_i - \gamma_1 - (X_i-\oX)^\top \beta_1\}^2  \text{ and }  \min_{\gamma_0, \beta_0} \sum_{i=1}^N (1-Z_i)\{Y_i - \gamma_0 - (X_i-\oX)^\top \beta_0\}^2.
\end{align}
It is not difficult to see that the least squares estimators for $\beta_1$ and $\beta_0$ from \eqref{eqn:two-regressions} are actually $\hat{\beta}_1$ and $\hat{\beta}_0$ defined before. 
Consequently, the resulting covariate-adjusted estimator $\htau(\hat{\beta}_1, \hat{\beta}_0)$ is equivalent to Lin's estimator $\htau_\textsc{l}$. In addition, the corresponding variance estimator constructed as in \eqref{eqn:hvN-beta} is asymptotically equivalent to the EHW variance estimator suggested by \citet{lin2013agnostic}.

From the above, Lin's estimator not only achieves the minimum true variance among all linearly covariate-adjusted estimators in \eqref{eqn:linear-estimators}, but also achieves the minimum estimated variance when we use the conservative variance estimator of form \eqref{eqn:hvN-beta}. 
A subtle issue here is that Lin's estimator uses estimated coefficients rather than fixed ones. 
With the technical tools discussed later, we can prove that the difference between Lin's estimator and the one with the oracle adjustment coefficients are asymptotically equivalent; see, e.g., \citet{li2017general}. 

% It can be shown that \eqref{eqn:two-regressions} is equivalent to the joint program \eqref{eqn:lin}.

% Another interpretation for Lin's estimator is based on adjusting for covariate imbalance. Let $\htau_X = \hX(1) - \hX(0)$ be the difference-in-means for the covariates. The linearly-adjusted estimator \eqref{eqn:linear-estimators} has an equivalent form
% \begin{align*}
%     \htau(\beta_1, \beta_0) = \htau - \gamma^\top \htau_X, \text{ where } \gamma = \frac{N_0}{N}\beta_1 + \frac{N_1}{N}\beta_0. 
% \end{align*}
% Similarly, Lin’s estimator has an equivalent form
% \begin{align*}
%     \htau_\scl = \htau - \hgamma^\top \htau_X, \text{ where } \hgamma = \frac{N_0}{N}\hbeta_1 + \frac{N_1}{N}\hbeta_0.
% \end{align*}
% We refer readers to \citet[Chapter 6]{ding2023first} for a more comprehensive discussion of different interpretations. 
% % \textbf{Some extensions.} 

\subsubsection{Further extensions}\label{sec:extensions}
There are a variety of extensions of covariate adjustment beyond treatment-control CREs.  

First, it is natural to consider generalization to multiple treatment levels. \citet{lu2016covariate} studied covariate adjustment in $2^K$ factorial designs by extending \eqref{eqn:two-regressions} to multi-level settings. \citet{zhao2023covariate} considered covariate adjustment in general multi-level experiments and made comprehensive comparison among the unadjusted estimator, Fisher's ANCOVA, and Lin's estimator. The unadjusted estimator is given by the regression:
\begin{align}\label{eqn:neyman-Q}
    (\hgamma_{1,\scn}, \ldots, \hgamma_{Q,\scn}) 
    = 
    \argmin_{\gamma_1, \dots, \gamma_Q} 
    \sum_{i=1}^N \lt(Y_i - \sum_{q=1}^Q\gamma_q \ind{Z_i = q}
    \rt)^2.
\end{align}

The generalization of Fisher's ANCOVA is given by the following \textit{additive treatment regression}: 
\begin{align}\label{eqn:fisher-Q}
    (\hgamma_{1,\scf}, \ldots, \hgamma_{Q,\scf},  \heta_{\scf}) = \argmin_{\gamma_1, \dots, \gamma_Q, \eta} \sum_{i=1}^N \lt\{Y_i - \sum_{q=1}^Q\gamma_q \ind{Z_i = q} - 
    (X_i-\oX)^\top \eta
    \rt\}^2.
\end{align}
Meanwhile, Lin's estimator can be generalized from either the regression with interaction perspective \eqref{eqn:lin} or the (estimated) variance minimization perspective \eqref{eqn:linear-estimators} or \eqref{eqn:two-regressions}. Here we present the former one, which applies the following fully interacted regressions: 
\begin{align}\label{eqn:lin-Q}
    (\hgamma_{1,\scl}, \ldots, \hgamma_{Q,\scl}, \heta_{1,\scl}, \ldots, \heta_{Q,\scl}) = \argmin_{\gamma_1, \dots, \gamma_Q, \eta_1, \dots, \eta_Q} \sum_{i=1}^N \lt\{Y_i - \sum_{q=1}^Q\gamma_q \ind{Z_i = q} - 
    \sum_{q=1}^Q \ind{Z_i = q} (X_i-\oX)^\top \eta_q
    \rt\}^2.
\end{align}

%\citet{ye2022toward} studied regression adjustment in multi-level experiments for arbitrary functions of response means (including linear contrasts, ratios, and odds ratios), multiple arms, guaranteed efficiency gain, optimality, and universal applicability
With the vectorized slope estimates $\hgamma_{*} = (\hgamma_{1,*}, \dots, \hgamma_{Q,*})^\top$, where $ * = \scn, \scf, \scl$, an estimator for the target average effect \eqref{eqn:ATE-Q} is given by the plug-in estimator
\begin{align}
    \htau_{*} = F^\top \hgamma_{*}, \quad * = \scn, \scf, \scl. 
\end{align}
Besides, we can obtain EHW variance estimators $\hV_{\mathrm{EHW}, *}$, which is conservative in large samples. In multi-level CRE, Lin's estimator is also guaranteed to be at least as efficient as Fisher's ANCOVA and Neyman's difference-in-means estimator.

Second, covariate adjustment has also been discussed in treatment-control trials when the dimension of the covariates is diverging or high-dimensional. For example, \citet{lei2020regression} proposed the following debiased estimator in treatment-control experiments:
\begin{gather}
    \htau_{\text{adj}}^{\text{de}} 
    = 
    \htau_\scl - 
    \lt(\frac{N_1}{N_0}\hDelta_0 - \frac{N_0}{N_1}\hDelta_1\rt).
\end{gather}
Here $\hDelta_z = N_z^{-1}\sum_{Z_i = z} \hat{e}_i H_{ii}$, $z= 0,1$, where $\hat{e}_i$ is the $i$-th residual from Lin's estimator \eqref{eqn:lin} and $H_{ii}$ is the $i$-th diagonal element of the hat matrix $H = X (X^\top X)^{-1}X^\top$, where $X$ is an $N \times p$ matrix with rows consisting of the covariates for the $N$ units. Under some structural conditions, the estimator $\htau_{\text{adj}}^{\text{de}}$ achieves asymptotic normality if the following condition holds:
\begin{align}
    \kappa^2 p \log p = o(1), \text{ where }\kappa = \max_{i=1,\dots,N} H_{ii}. 
\end{align}
In the favorable case where $\kappa = O(p/N)$, the dimension $p$ is allowed to grow as fast as $o(N^{2/3}/\log(N)^{1/3})$, which is a strictly weaker restriction than that of $\htau_\scl$. See also ref. \cite{lu2023debiased} for some recent development that allows $p$ to be in the same order as $N$.
As another example, \citet{bloniarz2016lasso} considered LASSO estimator for covariate adjustment in the high dimensional regime. Under a sparse linear model and some regularity conditions, the LASSO-adjusted regression estimator is asymptotically normal and the asymptotic variance is no greater than that of the difference-in-means estimator.

Third, some works explored the other variants of Lin's estimators. For example, \citet{zhao2023covariate} studied restricted least squares (RLS) and established for
the first time its properties for inferring ATE from the design-based perspective. \citet{guo2021generalized} considered generalized 
Oaxaca--Blinder
% oaxaca blinder 
estimators and extend the covariate adjustment framework from linear models to nonlinear ones; see also ref. \cite{cohen2020no}.

% \citet{ye2022toward}

\subsection{Rerandomization}\label{sec:rerandomization}

\citet{neyman1923application} focused on the CREs, 
which can balance all potential confounding factors, no matter observed or unobserved, on average and justifies the intuitive difference-in-means estimator for estimating the average treatment effect. 
In practice, in the design stage of an experiment, 
we often have access to a (rich) set of pretreatment covariates, 
and it has been a routine to check whether these covariates are balanced between different treatment groups. 
As commented by \citet{morgan2012rerandomization}, 
for a realized treatment allocation, 
the covariates are likely to be imbalanced; 
for example, with $10$ independent covariates,
% the probability that 
at least one of the $t$-statistics for checking the imbalance of these covariates will exceed $2$ with a probability of about $40\%$. 
It is then natural to incorporate the pretreatment covariate information into the design, aiming to get more balanced treatment groups as well as more efficient inference for treatment effects. 

Blocking is a classical and popular approach that can balance a few discrete covariates, but its implementation is not obvious once we have many continuous covariates. 
Rerandomization, a design recently formally proposed by \citet{morgan2012rerandomization}, provides a general approach to improve covariate imbalance, 
% solution to this issue of chance imbalance of covariates. 
although its idea has existed for a long time in the literature and dates back {to many earlier works \citep{sprott1993randomization, rubin2008comment, worrall2010evidence, cox2009randomization, bruhn2009pursuit, maclure2006measuring}.} 
In a recent survey of researchers conducting randomized experiments in developing countries \citep{Bruhn:2009}, the authors discovered that rerandomization has been commonly used in practice. 
For example, \citet{Lee2021} conducted a rerandomized experiment to study the effect of mobile banking for rural households and their migrated family members.

Under a general rerandomization design, 
for a randomly drawn treatment allocation, 
we will check the covariate balance between different treatment groups and see whether it satisfies a prespecified covariate balance criterion; 
if the balance criterion is met, we proceed to the actual experiment with this treatment allocation; otherwise, we redraw the treatment allocation and will keep redrawing until the balance criterion is met. 
% In order to facilitate covariate balance in the design stage, one strategy is called ``rerandomization''.
Although the balance criterion can be general, 
in the context of a treatment-control experiment, 
\citet{morgan2012rerandomization} suggested a balance criterion based on the Mahalanobis distance: 
\begin{align*}
    M & = 
    \hat{\tau}_X^\top \{ \Cov(\hat{\tau}_X) \}^{-1} 
    \hat{\tau}_X
    = 
    \frac{N_1 N_0}{N} 
    \hat{\tau}_X^\top 
    (S_X^2)^{-1} \hat{\tau}_X, 
\end{align*}
recalling that $N_1$ and $N_0$ are the treated and control group sizes, 
$\hat{\tau}_X$ is the difference-in-means of covariates defined as in Section \ref{sec:Lin}, 
and $S_X^2$ is the finite population covariance matrix of covariates defined as follows:
\begin{align*}
    % \htau_X = \hX(1) - \hX(0), \quad 
    S_X^2 = \frac{1}{N-1}\sum_{i=1}^N (X_i - \oX)(X_i - \oX)^\top. 
\end{align*}
Under rerandomization using the Mahalanobis distance, denoted by ReM, 
we will repeatedly draw random treatment assignment from the CRE until getting an acceptable one with the corresponding Mahalanobis distance bounded by a prespecified threshold $a$. 

Importantly, the analysis for rerandomization needs to take into account the selection step in its design. 
This is often ignored in practice, and rerandomization is often analyzed as if it was a CRE. 
\citet{morgan2012rerandomization} proposed randomization tests for sharp null hypotheses, 
employing assignments generated randomly in accordance with the rerandomization protocol.
% with assignments randomly generated following the rerandomization protocol. 
More recently, \citet{li2018asymptotic} conducted Neyman-type large-sample inference for rerandomization, considering also the intuitive difference-in-means estimator. 
They demonstrated that, asymptotically, the difference-in-means estimator is more concentrated around the true average treatment effect with smaller asymptotic variance and shorter asymptotic quantile ranges, 
and proposed accurate confidence intervals for the average effect, which are always shorter than Neyman's intervals for the CRE while remaining valid asymptotically under ReM.

In recent years, rerandomization has been extended to more general experiments, such as factorial experiments \citep{branson2016improving, li2020Arerandomization}, blocked experiments \citep{wang2023rerandomization, johansson2022rerandomization}, 
and survey experiments \citep{YQL2021}, 
and it can also be combined with covariate adjustment discussed in Section \ref{sec:cov-adj} \citep{li2020Brerandomization}. \citet{zhao2024no} studied the procedure of conducting rerandomization directly based on
p-values from covariate balance tests, which is a general strategy that works for many basic designs. 
An alternative rerandomization scheme that randomizes treatment assignments multiple times and chooses the one with the best covariate balance has also been used in practice \citep{Bruhn:2009}, and its property has recently been studied in ref. \cite{wang2023asymptotic}.

% Moreover, regression adjustment can be performed to further improve efficiency, as shown by \citet{li2020Brerandomization}. The discussion is also extended to factorial experiments \citep{li2020Arerandomization}.

% \begin{table}[!ht]
% \centering
% \caption{Design and analysis of experiments \citep{ding2023first}}\label{tb::design-analysis-2x2}
% \label{tb::rerandomization-regression}
% \begin{tabular}{c|cccc}
%   & & & analysis & \\
%   \hline 
%   & CRE & $\hat\tau$ \citep{neyman1923application} &  $\stackrel{1}{\longrightarrow}$  &  $\hat\tau_\textsc{L}$ \citep{lin2013agnostic} \\
%   design&  &  $2\Big\downarrow$ && $  \Big\downarrow 4$ \\
%   & ReM & $\hat\tau$  \citep{li2018asymptotic} &  $\stackrel{3}{\longrightarrow}$  &$\hat\tau_\textsc{L}$ \citep{li2020Brerandomization} \\
% \end{tabular}
% \end{table}

% {
% \color{red} For Xinran to write
% }

% \citet{morgan2012rerandomization},
% \citet{li2018asymptotic}, \citet{li2020Arerandomization}, \citet{li2020Brerandomization}, 
% \citet{wang2022rerandomization}

% \section{Permutational central limit theorems and Berry-Esseen bounds with application in randomization inference}\label{sec:valid-inference}

\section{Permutational central limit theorems}\label{sec:PCLT}

With all the design and analysis strategies introduced above, one natural question is how to theoretically justify their statistical property. In the following two sections, we focus on the technical aspect of CREs. The main question to answer is how to deliver valid inference with different estimators for different designs. Permutational/combinatorial CLTs and BEBs are core to the technical development of randomization-based inference.  in Sections \ref{sec:PCLT} and \ref{sec:PBEB}, we summarize the theoretical results regarding permutational CLTs and BEBs and discuss their application in analyzing randomized experiments.

\subsection{Sample average under simple random sampling}
We start with the simple random sampling from a finite population \citep{cochran1977sampling}. Let $\{a_N(i)\}_{i=1}^N$ be a sequence of real numbers. Suppose we randomly sample $N_1$ elements \textit{without replacement} from the population and use a binary variable $Z_i$ to indicate the sampling status of the $i$-th element, i.e., $Z_i = 1$ indicates $a_N(i)$ being sampled while $Z_i = 0$ not sampled. Write $N_0 = N - N_1$. Consider the sample average obtained from the above procedure:
\begin{align}\label{eqn:srs}
  \Gamma = \frac{1}{N_1}\sum_{i=1}^N a_N(i) \ind{Z_i = 1}.
\end{align}
$\Gamma$ has mean and variance%\todo{is $=$ used before?}
\begin{align}
    \E{\Gamma} = \oa_N, \quad   V_N = \Var{\Gamma} = (\frac{1}{N_1} - \frac{1}{N}) S^2_{N},
\end{align}
where
\begin{align}
    \oa_N = \frac{1}{N}\sum_{i=1}^N a_N(i), \quad S^2_{N} =  \frac{1}{N-1}\sum_{i=1}^N (a_N(i) - \oa_N)^2.
\end{align}

One fundamental technical question is to establish central limit theorems for $\Gamma$ to characterize its asymptotic distribution. \citet{erdos1959central} established the following CLT for \eqref{eqn:srs}:
\begin{proposition}\label{prop:srs-clt}
    If for any $\epsilon > 0$,
    \begin{align}\label{eqn:lf}
        \lim_{N \to \infty} \frac{\sum_{i=1}^N (a_N(i) - \oa_N)^2\ind{|a_N(i) - \oa_N| > \epsilon N_1\sqrt{V_N}}}{\sum_{i=1}^N (a_N(i) - \oa_N)^2} \to 0,
    \end{align}
    then as $N\to\infty$,
    \begin{align}
        \frac{\Gamma - \E{\Gamma}}{\sqrt{V_N}} \rightsquigarrow \cN(0,1).
    \end{align}
\end{proposition}

\citet{hajek1960limiting} further proved that Condition \eqref{eqn:lf} is not only sufficient but also necessary provided that $N_1,N_0\to\infty$.  Moreover, Theorem \ref{prop:srs-clt} covers some other works on finite population CLTs. For example, \citet{madow1948limiting} proved asymptotic normality under the conditions that $N_1\to\infty$ and 
there exists $\delta\in(0,1)$ such that $N_1/N < 1 - \delta$ when $N$ is sufficiently large, and that the following moment condition holds:
\begin{align}\label{eqn:madow-moment}
    \frac{N^{-1}\sum_{i=1}^N |a_N(i) - \oa_N|^r}{\{N^{-1}\sum_{i=1}^N (a_N(i) - \oa_N)^2\}^{r/2}}  = O(1), \text{ for all integers } r > 2.
\end{align}
The above moment condition \eqref{eqn:madow-moment} is stronger than \eqref{eqn:lf}, because for any $r>2$, 
\begin{align}
    &\frac{\sum_{i=1}^N (a_N(i) - \oa_N)^2\ind{|a_N(i) - \oa_N| > \epsilon N_1\sqrt{V_N}}}{\sum_{i=1}^N (a_N(i) - \oa_N)^2}\\
    =&
    {(1 - \frac{N_1}{N})\frac{N_1 \epsilon^2 }{N-1}\sum_{i=1}^N (\frac{a_N(i) - \oa_N}{\epsilon N_1\sqrt{V_N}})^2\ind{|a_N(i) - \oa_N| > \epsilon N_1\sqrt{V_N}}} \\
    \le   &
    {(1 - \frac{N_1}{N})\frac{N_1 \epsilon^{2-r} }{N-1}
    \sum_{i=1}^N \lt|\frac{a_N(i) - \oa_N}{ N_1\sqrt{V_N}}\rt|^r} \\
    \lesssim &
    \frac{1}{N_1^{\frac{r}{2} - 1}}\cdot\frac{N^{-1}\sum_{i=1}^N |a_N(i) - \oa_N|^r}{\{N^{-1}\sum_{i=1}^N (a_N(i) - \oa_N)^2\}^{r/2}},
\end{align}
which converges to zero under \eqref{eqn:madow-moment}.

\citet{david1938limiting}  established a CLT for the hypergeometric distribution, which is a special case of \citet{madow1948limiting} thus also stronger than \eqref{eqn:lf}. \citet[Section 2.1]{li2017general} also provided a thorough exposition of CLT under the simple random sampling scheme with a sufficient condition based on the maximum squared distance. 
% The work by \citet{sen1970hajek} established the H\'{a}jek-R\'{e}nyi inequality, which can be used to quantify the almost sure convergence rate of the sampling statistics.

\subsection{Simple linear rank statistics}

The sample average in \eqref{eqn:srs} from a SRS is a special case of a more general type of permutational statistics, called \textit{simple linear rank statistics}.  Formally,
let $\{a_N(i)\}_{i=1}^N$ and $\{b_N(i)\}_{i=1}^N$ be two sequences of real numbers. Let $\pi$ be a random permutation over the indices $1,\dots, N$, with $\pi(i)$ denoting the permuted index of $i$. A simple linear rank statistic is defined as
\begin{align}\label{eqn:slrs}
    \Gamma = \sum_{i=1}^N a_N(i) b_N(\pi(i)),
\end{align}
which has mean and variance
\begin{align}
    \E{\Gamma} = N \oa_N\cdot \ob_N, \quad V_N = \Var{\Gamma} =  \frac{1}{N-1} \cdot \sum_{i=1}^N(a_N(i) - \oa_N)^2\cdot \sum_{i=1}^N(b_N(i) - \ob_N)^2.
\end{align}
In particular, if we take $b_N(i) = N_1^{-1}$ for $i = 1,\dots, N_1$ and  $b_N(i) = 0$ for $i = N_1 + 1,\dots, N$, then \eqref{eqn:slrs} gives the sample average \eqref{eqn:srs} in SRS.
The statistic in 
\eqref{eqn:slrs} has been studied by many researchers. \citet{wald1944statistical} established CLT under the following condition: for all integers $r > 2$,
% define the quantity
% Proposition (Theorem), present the weakest version
\begin{align}\label{eqn:slrs-ww}
    \frac{N^{-1}\sum_{i=1}^N (a_N(i) - \oa_N)^r}{\{N^{-1}\sum_{i=1}^N (a_N(i) - \oa_N)^2\}^{r/2}}  = O(1), \quad \frac{N^{-1}\sum_{i=1}^N (b_N(i) - \ob_N)^r}{\{N^{-1}\sum_{i=1}^N (b_N(i) - \ob_N)^2\}^{r/2}}  = O(1).
\end{align}
\citet{noether1949theorem} proved CLT under the following condition that is slightly weaker than \citet{wald1944statistical}: for all integers $r > 2$, 
\begin{align}\label{eqn:slrs-noether}
    \frac{N^{-1}\sum_{i=1}^N (a_N(i) - \oa_N)^r}{\{N^{-1}\sum_{i=1}^N (a_N(i) - \oa_N)^2\}^{r/2}}  = O(1), \quad \frac{\sum_{i=1}^N (b_N(i) - \ob_N)^r}{\{\sum_{i=1}^N (b_N(i) - \ob_N)^2\}^{r/2}}  = o(1),
\end{align}
% \eqref{eqn:slrs-noether} 
which, however, 
is not symmetric for $a_N(i)$'s and $b_N(i)$'s. 
\citet{hoeffding1951combinatorial} further proved CLT under a weaker and symmetric condition: for all integers $r>2$,
\begin{align}\label{eqn:slrs-hoeffding}
    N^{\frac{r}{2} - 1}\frac{\sum_{i=1}^N (a_N(i) - \oa_N)^r}{\{\sum_{i=1}^N (a_N(i) - \oa_N)^2\}^{r/2}} 
\cdot \frac{\sum_{i=1}^N (b_N(i) - \ob_N)^r}{\{\sum_{i=1}^N (b_N(i) - \ob_N)^2\}^{r/2}}  = o(1).
\end{align}
\citet{motoo1956hoeffding} proved that CLT holds under an even weaker Lindeberg type condition: 
\begin{proposition}\label{prop:srls-motoo}
    Suppose that for any $\epsilon > 0$, 
\begin{align}\label{eqn:srls-motoo}
\lim_{N\to\infty}\frac{\sum_{i,j=1}^N (a_N(i) - \oa_N)^2(b_N(j) - \ob_N)^2\ind{|(a_N(i) - \oa_N)(b_N(j) - \ob_N)| > \epsilon \sqrt{V_N} }}{\sum_{i,j=1}^N (a_N(i) - \oa_N)^2(b_N(j) - \ob_N)^2} = 0. 
\end{align}
Then 
\begin{align*}
    \frac{\Gamma - \E{\Gamma}}{\sqrt{\Var{\Gamma}}} \rightsquigarrow \cN(0,1). 
\end{align*}
\end{proposition}

\citet{hajek1961some} proved further that Condition \eqref{eqn:srls-motoo} is not only sufficient but also necessary, 
% . Also \citet{hajek1961some} 
and 
presented a comprehensive comparison of the CLT conditions introduced in the literature. There are also several multidimensional extensions based on the Cram\'{e}r-Wold device. See, for example, refs \cite{fraser1956vector}, \cite[Section 7]{hajek1961some}, 
and \cite[Lemma S.3.3]{diciccio2017robust}. 
% , etc. 
%\citet{dwass1953asymptotic} 

\subsection{General univariate linear permutational statistics}\label{sec:linear_perm_clt}
% \textbf{Univariate case.}
Taking one step further from the simple linear rank statistics, the permutational central limit theorems are 
% proved 
proposed
for the following {\it linear permutational statistic}:
\begin{align}\label{eqn:sips}
    \Gamma = \sum_{i=1}^N M_N(i,\pi(i)),
\end{align}
where $\{M_N(i,j)\}_{i,j\in[N]}$ is a matrix in $\bbR^{N\times N}$. In particular, if we take $M_N(i,j) = a_N(i)b_N(j)$, \eqref{eqn:sips} recovers \eqref{eqn:slrs}.   \citet{hoeffding1951combinatorial} computed the mean and variance of \eqref{eqn:sips}:
\begin{align*}
    \E{\Gamma} = \frac{1}{N}\sum_{i,j=1}^N M_N(i,j),\quad V_N = \Var{\Gamma} = \frac{1}{N-1}\sum_{i,j=1}^N \tilde{M}_N(i,j)^2.
\end{align*}
Here $\tilde{M}_N(i,j)$ is the centered array based on the following rule:
\begin{align}\label{eqn:tildeM}
    \tilde{M}_N(i,j) = M_N(i,j) - \frac{1}{N}M_N(i,+) - \frac{1}{N}M_N(+,j) + \frac{1}{N^2} M_N(+,+).
\end{align}
where ``$+$'' means summation over the corresponding index.

Moreover, \citet{hoeffding1951combinatorial} 
% was able to show 
showed that the asymptotic normality of $\Gamma$ in \eqref{eqn:sips} holds under the following condition: 
\begin{align}\label{eqn:sips-hoeffding}
\lim_{N\to\infty}\frac{N^{-1}\sum_{i,j=1}^N \tilde{M}_N(i,j)^r}{\{N^{-1}\sum_{i,j=1}^N \tilde{M}_N(i,j)^2\}^{r/2}} = 0, \quad \text{ for all integers $r >2$}. 
\end{align}
Condition \eqref{eqn:sips-hoeffding}  is equivalent to Condition \eqref{eqn:slrs-hoeffding} in the simple linear rank statistics setting. 
A more compact sufficient condition for \eqref{eqn:sips-hoeffding} is also provided in \citet{hoeffding1951combinatorial}:
\begin{align}\label{eqn:suff-sips-hoeffding}
\lim_{N\to\infty}\frac{\max_{i,j\in[N]} \tilde{M}_N(i,j)^2}{N^{-1}\sum_{i,j=1}^N \tilde{M}_N(i,j)^2} = 0.
\end{align}
\citet{motoo1956hoeffding} weakened \citet{hoeffding1951combinatorial}'s condition \eqref{eqn:sips-hoeffding} to the following Lindeberg type condition:
\begin{proposition}[Main theorem of \citet{motoo1956hoeffding}]\label{prop:clt-motoo}
    Suppose for any $\epsilon > 0$, 
    \begin{align}
      \lim_{N\to\infty} \frac{\sum_{i,j=1}^N \tilde{M}_N(i,j)^2\ind{|\tilde{M}_N(i,j)| > \epsilon  \sqrt{V_N}}}{\sum_{i,j=1}^N \tilde{M}_N(i,j)^2} = 0.
    \end{align}
    Then 
    \begin{align*}
       \frac{\Gamma - \E{\Gamma}}{\sqrt{\Var{\Gamma}}} \rightsquigarrow \cN(0,1). 
    \end{align*}
\end{proposition}
\begin{remark}\label{rmk:compare-PCLTs}
    Although Proposition \ref{prop:clt-motoo} gives the weakest condition for permutational CLT in the literature, it is not very convenient for use in many concrete examples. On the contrary, Condition \eqref{eqn:suff-sips-hoeffding} involves the maximum of the centered matrices and is simpler for use and interpretation. Condition \eqref{eqn:suff-sips-hoeffding} and its multivariate generalization (presented in  \eqref{eqn:suff-mv-sips-fraser} in Section \ref{sec:mv-sips}) are frequently applied to investigation of the properties of various analysis and design strategies in randomized experiments presented in Section \ref{sec:design-analysis}; 
    {\lxr for example, they can faciliate the proof for the convergence of variance estimation.} 
    We will have more discussion in Section \ref{sec:apply-PCLT}.  
\end{remark}

\subsection{General multivariate linear permutational statistics}\label{sec:mv-sips}

% \textbf{Multivariate case.} 

We now discuss the generalization of \eqref{eqn:sips-hoeffding} to a multivariate case. Concretely, define the multivariate linear permutational statistics:
\begin{align}\label{eqn:mv-sips}
    \Gamma = (\Gamma_1, \dots, \Gamma_H)^\top, \quad \Gamma_h =  \sum_{i=1}^N M_{N,h}(i,\pi(i)),
\end{align}
where $\{M_{N,h}(i,j)\}_{i,j\in[N]}$, $h=1,\dots,H$ are $H$ matrices in $\bbR^{N\times N}$. \citet[Appendix A.1.]{shi2022berry} presented many basic facts about \eqref{eqn:mv-sips}, including its mean and covariance calculation and its standardization. 
\citet{fraser1956vector} extended \citet{hoeffding1951combinatorial} to the multi-dimensional setting by applying the Cram\'{e}r-Wold device to establish a multivariate central limit theorem. More concretely, define the centered version of $\tilde{M}_{N,h}$ in the same way as \eqref{eqn:tildeM}. %\todo{Just to confirm, this is different from the standardization later?} 
\citet{fraser1956vector} proposed the following condition for CLT as an extension to \eqref{eqn:sips-hoeffding}:
\begin{align}\label{eqn:mv-sips-fraser}
\lim_{N\to\infty}\frac{N^{-1}\sum_{i,j=1}^N \tilde{M}_{N,h}(i,j)^r}{\{N^{-1}\sum_{i,j=1}^N \tilde{M}_{N,h}(i,j)^2\}^{r/2}} = 0, \quad \text{ for all integers $r >2$ and $h \in [H]$}. 
\end{align}
Similarly, \citet{fraser1956vector} also provided a sufficient condition for \eqref{eqn:mv-sips-fraser}:
\begin{align}\label{eqn:suff-mv-sips-fraser}
    \lim_{N\to\infty}\frac{\max_{i,j\in[N]} \tilde{M}_{N,h}(i,j)^2}{N^{-1}\sum_{i,j=1}^N \tilde{M}_{N,h}(i,j)^2} = 0, \quad \text{ for } h\in[H]. 
\end{align}
The condition in \eqref{eqn:suff-mv-sips-fraser} is further utilized 
% and extended in 
by 
\citet{li2017general} to build asymptotic normality results for 
% many permutational related problems. 
{\lxr analyzing treatment effects in multi-level CREs}.

\subsection{Application of permutational CLT in randomization-based inference}\label{sec:apply-PCLT}

In this subsection, we collect several theoretical arguments in randomization-based inference that apply permutational CLTs to deliver technical justification for studying average treatment effects. 

% \textbf{Wald-type inference in CRE}. 

\subsubsection{Wald-type inference in completely randomized experiments}\label{sec:inference-CRE}

{\lxr Consider analyzing a multi-level CRE discussed in Section \ref{sec:multi_CRE}, and adopt the notation introduced there.}
% For the confidence region \eqref{eqn:CR}, 
\citet{li2017general} proved the following result to justify 
% its 
the asymptotic validity 
of the  confidence region \eqref{eqn:CR}
under several regularity conditions. 
\begin{proposition}[Theorem 5 and Proposition 3 of \citet{li2017general}]\label{prop:wald-inference-CRE}
Let $Q$ be fixed and $N$ go to infinity. If the covariance matrix $S$ have limiting values, $N_q/N$ has positive limiting value, and 
$$
\max_{1\le q\le Q} \max_{1\le i\le N} |Y_i(q) - \oY(q)|^2/N \to 0,
$$ 
then the following conclusions hold:
\begin{enumerate}[label = \em(\roman*)]
    \item {\em Asymptotic normality.} $N\Var{\htau}$ has a semi-positive definite limiting value $V_\infty$, and 
\begin{align}
    \sqrt{N}(\htau - \tau) \rightsquigarrow \cN(0, V_\infty),
\end{align}
where $\tau$ and $\hat{\tau}$ are the average treatment effect and the corresponding estimator in \eqref{eqn:ATE-Q} and \eqref{eqn:htau-Q}.
    \item {\em Variance estimation.} The 
    % empirical variance estimator
    sample variance 
    $\hS(q,q)$ in \eqref{eq:S_hat_qq} is consistent for $S(q,q)$ in \eqref{eq:S_qq}. 
    
    \item {\em Wald-type inference.} If the limit of $ N F^\top \diag{N_q^{-1} S(q,q)} F $ is nonsingular, then $\hV$ in \eqref{eqn:hVN} is nonsingular with probability converging to one, and the Wald-type confidence region \eqref{eqn:CR} has asymptotic coverage rate at least $1-\alpha$. Moreover, the asymptotic coverage rate equals $1 - \alpha$ if and only if the causal effects are asymptotically additive, in the sense that $\lim_{N\to\infty} F^\top S F = 0$.
\end{enumerate}

\end{proposition}

We briefly comment on the technical details behind Proposition \ref{prop:wald-inference-CRE}. Proposition \ref{prop:wald-inference-CRE}(i) utilized the permutational CLT for general linear permutational statistics. In particular, the estimate $\htau$ 
% can be written in the form of 
follows the same distribution as 
\eqref{eqn:mv-sips} with matrices $M_{N,h}$ defined as follows: 
% let $N_0 = 0$, 
for $i,j \in [N]$,
\begin{align}
    M_{N,h}(i,j)
    = N_q^{-1} F(q, h) Y_i(q), \quad N_{q-1} + 1 \le j \le N_q,
\end{align}
where $N_0 = 0$ and $F(q,h)$ is the $(q,h)$th element of the contrast matrix $F$. 
Now applying Condition \eqref{eqn:suff-mv-sips-fraser}, we can obtain the regularity conditions on potential outcomes in Proposition \ref{prop:wald-inference-CRE} and justify the multivariate asymptotic normality of the estimator $\htau$. Proposition \ref{prop:wald-inference-CRE}(ii) applied the Chebyshev inequality to the sample variance estimators and showed their consistency. Proposition \ref{prop:wald-inference-CRE}(iii) combined (i) and (ii) and formally established the asymptotic validity of \eqref{eqn:CR}. As a side remark, the matrix $V_\infty$ is not required to be invertible, because the multivariate combinatorial CLT proved by ref. \cite{fraser1956vector} does not require an invertible limit for the covariance matrix. However, to justify the validity of the Wald-type confidence intervals in part (iii), invertibility is required. 

Proposition \ref{prop:wald-inference-CRE} covers the treatment-control experiments as special cases.  In other words, under certain regularity conditions, the difference-in-means estimator $\htau$ in \eqref{eqn:htau} is asymptotically normal, and 
the variance estimator $\hat{V}$ in \eqref{eqn:vN} is consistent for a limit that is no less than the true asymptotic variance of $\hat{\tau}$.
% and $\hv$ is consistent.
% \begin{align}\label{eqn:asymp-htauN-hvN}
%     \frac{\htau - \tau}{\sqrt{\Var{\htau}}} \rightsquigarrow \cN(0,1), \quad \hv - \E{\hv} \xrightarrow{\bbP} 0.
% \end{align}
These then justify the asymptotic validity of the level-$(1-\alpha)$ confidence interval in \eqref{eqn:CI}.

\subsubsection{Analyzing covariate adjustment}

For covariate adjustment, we present a result by \citet{zhao2023covariate}. Let $\hY_* \in \mathbb{R}^Q$ ($* = \scn, \scf, \scl$) be the estimators for the averaged potential outcomes across all $Q$ treatment levels from three estimation strategies: $\hY_\scn$ from Neyman's approach, $\hY_\scf$ from Fisher's ANCOVA and $\hY_\scl$ from Lin's regression. These three estimators correspond to the coefficients $\hgamma_*$ in front of the treatment indicators from the regressions introduced in Section \ref{sec:extensions}, i.e., \eqref{eqn:neyman-Q}, \eqref{eqn:fisher-Q} and \eqref{eqn:lin-Q}. 
We slightly modified the notation  
in order to better present the results. 
% next Proposition, 
Also, define $\hat{V}_*$ to be the corresponding EHW robust covariance estimator from the three regressions. The following proposition 
% \ref{prop:asymp-X-adj} 
from ref. \cite{zhao2023covariate} established the asymptotic properties of these point and variance estimators.
\begin{proposition}[Lemma 1 of \citet{zhao2023covariate}]\label{prop:asymp-X-adj}
    Let $N \to \infty$. 
    % Assume the covariates are centered: $N^{-1}\sum_{i=1}^N X_i = 0$. 
    Assume that, for $q\in[Q]$, $e_q = N_q/N$ has a limit in $(0, 1)$. Assume
    that the first two finite population moments of $\{Y_i(q), X_i, X_iY_i(q): q \in [Q] \}$ have finite limits, and both $S^2_X = (N-1)^{-1}\sum_{i=1}^N X_i X_i^\top$ and its limit are nonsingular, where the covariates have been centered so that $N^{-1}\sum_{i=1}^N X_i = 0$. Also assume that $N^{-1}\sum_{i=1}^N Y_i^4(q) = O(1) $, $ N^{-1}\sum_{i=1}^N \|X_i\|_4^4 = O(1) $, and $ N^{-1}\sum_{i=1}^N \|X_iY_i(q)\|_4^4 = O(1) $. Then the following results hold:
    \begin{enumerate}[label = \em (\roman*)]
        \item {\em Asymptotic normality}. $\sqrt{N}(\hY_* - \oY) \rightsquigarrow \cN(0, V_*) $ for some {\em $V_* \succeq 0$, $* = \scn, \scf, \scl$}. 

        \item {\em Conservative variance estimation}. {\em $ \mathrm{plim}_{N\to\infty} ~ N \hV_{*, \text{EHW}}  \succeq  V_*$, $* = \scn, \scf, \scl$.}
        
        \item {\em Efficiency comparison}. {\em  $ V_\scl \preceq V_\scn $ and $ V_\scl \preceq V_\scf $.} 
    \end{enumerate}
\end{proposition}

Proposition \ref{prop:asymp-X-adj}(i) established the asymptotic normality property of Neyman's difference-in-means, Fisher's ANCOVA and Lin's estimator. Together with the conservative variance estimation 
% conservativeness suggested by 
in Proposition \ref{prop:asymp-X-adj}(ii), one can justify the asymptotic validity of Wald-type confidence regions constructed from 
% covariate-adjusted 
these 
estimators. Proposition \ref{prop:asymp-X-adj}(iii) indicates that Lin's estimator guarantees at least as much asymptotic efficiency as the difference-in-means estimator and Fisher's ANCOVA.  

In terms of the technical derivation, Proposition \ref{prop:asymp-X-adj}(i) utilized Proposition \ref{prop:wald-inference-CRE}, which is a result motivated by Hoeffding \& Fraser's permutational CLT 
% under 
(more specifically, Conditions \eqref{eqn:suff-sips-hoeffding} and \eqref{eqn:suff-mv-sips-fraser}) and can accommodate vector outcomes and multi-level randomized experiments. In particular, if we study the pseudo potential outcome vector $ (Y_i(q), X_i^\top)^\top$, for  $q\in[Q]$, 
we can apply similar tricks as in Section \ref{sec:inference-CRE} to establish a multivariate CLT for the arm-wise sample means
$(\hY_q, \hX_q^\top)^\top $ for $q\in[Q] $.
based on Condition \eqref{eqn:suff-mv-sips-fraser}. 
The covariate-adjusted estimator $\hgamma_{\scl}$ in \eqref{eqn:lin-Q} can be formulated as linear combinations of these sample means, where the 
% weights 
combination coefficients are consistent for certain constant coefficients in the sense that their difference is of order $o_{\bbP}(1)$.
% converge in probability. 
Then a CLT can be derived after filling in the details \citep{zhao2023covariate}. Proposition \ref{prop:asymp-X-adj}(ii) utilized the Chebyshev inequality 
% to show consistency in probability
under the bounded moment conditions.  Proposition \ref{prop:asymp-X-adj}(iii) involves some delicate analysis of the limiting variance structure $V_{*}$ for $* = \scn, \scf, \scl$, which has closed-form 
% representations 
expressions \citep{zhao2023covariate}. 

\subsubsection{Analyzing rerandomization}\label{sec:rem_fix_a}

For rerandomization using the Mahalanobis distance discussed in Section \ref{sec:rerandomization}, 
we adopt the notation introduced there and present the following result by \citet{li2018asymptotic}. Define $V$ as the variance of $\sqrt{N}\hat{\tau}$ under the CRE, 
and $R^2$ as the squared multiple correlation between the difference-in-means of outcome and covariates; see Proposition 1 of ref. \cite{li2018asymptotic} for its explicit expression.
Let $\varepsilon_0\sim \mathcal{N}(0,1)$, $L_{K,a}\sim D_1 \mid D'D\le a$ with $D=(D_1, \ldots, D_K)^\top \sim \mathcal{N}(0, I_K)$, and $\varepsilon_0 \indep L_{K,a}$

\begin{proposition}[\citet{li2018asymptotic}, Theorems 1 and 2 and Appendix A4.2]\label{prop:rem}
Consider ReM with a fixed positive threshold $a$,
% and a fixed number of covariates $K$ that do not vary with the sample size, 
and assume that,  
as $N\rightarrow \infty$, 
(a) the proportions of units under treatment and control have positive limits, 
(b)
the finite population variances and covariances for potential outcomes and covariates have limits, and 
the limit of $S^2_X$ is nonsingular; 
and (c) $N^{-1} \max_{1\le i \le N} |Y_i(z) - \bar{Y}(z)|^2 \rightarrow 0$ and 
$N^{-1} \max_{1\le i \le N} \|X_i - \bar{X}\|_2^2 \rightarrow 0$. 
\begin{enumerate}[label = \em (\roman*)]
        \item {\em Asymptotic distribution}. 
        $\sqrt{N}(\hat{\tau} - \tau) \mid M\le a
        \ \rightsquigarrow \ 
        V^{1/2}(\sqrt{1-R^2}\cdot \varepsilon_0 + \sqrt{R^2}\cdot L_{K,a}).
        $
        % where the conditioning on the event that the Mahalanobis distance $M$ is bounded by $a$ is used to emphasize the ReM design,  
        % $V$ is the variance of $\sqrt{N}\hat{\tau}$ under the CRE, 
        % $R^2$ is an $R^2$-type measure for the association between potential outcomes and covariates, 
        % $\varepsilon_0\sim \mathcal{N}(0,1)$, 
        % $L_{K,a}\sim D_1 \mid D'D\le a$ with $D=(D_1, \ldots, D_K)^\top \sim \mathcal{N}(0, I_K)$, 
        % and $\varepsilon_0 \indep L_{K,a}$. %\todo{maybe simplify, just say see detailed definition in [REF]} 
      % {\lxr define $V$, }  
        \item {\em Conservative inference}. 
        We can construct estimators $\hat{V}$ and $\hat{R}^2$ such that,
        as $N\rightarrow \infty$, 
        $\mathrm{plim}_{N\to\infty}( \hat{V} - V ) \ge 0$ %\todo{see the definition of $\hat{V}$ and $\hat{R}^2$ in [REF]?} 
        and $\mathrm{plim}_{N\to\infty} ( \hat{V}\hat{R}^2 - VR^2 ) = 0$. 
        
        \item {\em Efficiency comparison}. 
        % The asymptotic distribution under ReM, $V^{1/2}(\sqrt{1-R^2}\cdot \varepsilon_0 + \sqrt{R^2}\cdot L_{K,a})$,   has smaller (or the same) variance and shorter (or the same) symmetric quantile ranges than that under the CRE, 
        The asymptotic distribution under ReM has a smaller (or equal) variance and narrower (or equal) symmetric quantile ranges than that under the CRE. Here, the symmetric quantile range means the interval formulated by the lower and upper $\alpha/2$-th quantile of the asymptotic distribution of $\sqrt{N}(\htau-\tau)$, for $\alpha \in (0,1)$.
    \end{enumerate}
\end{proposition}

Proposition \ref{prop:rem}(i) means that the distribution of the difference-in-means estimator under rerandomization converges in distribution to the same limit as $V^{1/2}(\sqrt{1-R^2}\cdot \varepsilon_0 + \sqrt{R^2}\cdot L_{K,a})$, a convolution of a Gaussian and a constrained Gaussian random variables, 
{\lxr where the coefficient depends crucially on $R^2$ that represents an $R^2$-type measure for the association between potential outcomes and covariates;}
% Moreover, $R^2$ represents an $R^2$-type measure for the association between potential outcomes and covariates; 
see ref. \cite{li2018asymptotic} for more details. 
% In Proposition \ref{prop:rem}, 
% (i) shows the asymptotic distribution of the difference-in-means estimator under rerandomization. 
Interestingly, unlike that under the CRE,  
the asymptotic distribution of $\hat{\tau}$ under rerandomization is non-Gaussian in general, 
while it is still symmetric and unimodal around zero \citep{li2018asymptotic}. 
In addition, when $a = \infty$ or $R^2=0$, 
the asymptotic distribution reduces to that for the CRE. 
The former is not surprising, because ReM without rejecting any assignment is essentially the CRE. 
The latter is also intuitive, 
implying that ReM using covariates that are irrelevant to potential outcomes is asymptotically equivalent to the CRE without using any covariates. 
In Proposition \ref{prop:rem}(ii), we omit the explicit expressions of the estimators for conciseness. 
As a side remark, \citet{li2018asymptotic} uses an estimator for $V$ that is less conservative than or asymptotically equivalent to \eqref{eqn:vN} by utilizing the covariate information; we refer interesting readers to \citep{li2018asymptotic,DFM2019} for details.
Importantly,
Proposition \ref{prop:rem}(ii) shows that we can consistently estimate the coefficient of $L_{K,a}$ in the asymptotic distribution, while only conservatively estimate the coefficient of $\varepsilon_0$. 
Fortunately, due to the symmetric and unimodal property of the asymptotic distribution, 
these will lead to conservative variance estimation and confidence intervals.

Proposition \ref{prop:rem}(iii) demonstrates the advantage of rerandomization over the CRE. 
In particular, the stronger the association between covariates and potential outcomes, as measured by $R^2$, the larger the gain from ReM \citep{li2018asymptotic}. 
\citet{Zach2023} recently extended the comparison to non-symmetric quantile ranges.

We now discuss the technical aspects of Proposition \ref{prop:rem}. 
A key for its derivation is to notice that the distribution of $\hat{\tau}$ under rerandomization is the same as its conditional distribution under the CRE given that the covariate balance criterion is satisfied. 
This is emphasized by the conditioning on $M\le a$ in Proposition \ref{prop:rem}(i). 
Thus, to understand this conditional distribution, 
it suffices to study the joint distribution of the difference-in-means vector $(\hat{\tau}, \hat{\tau}_X^\top)$ for both outcome and covariates, noting that $M$ is a deterministic function of $\hat{\tau}_X$. 
Such a joint distribution will be asymptotically normal, 
which can be derived using  Proposition \ref{prop:wald-inference-CRE} 
by viewing $(Y(z), X)$ as a ``potential outcome'' vector. 
For the asymptotic conservative  inference in  Proposition \ref{prop:rem}(ii), 
we can study 
the probability limits of the estimators $\hat{V}$ and $\hat{R}^2$ again  
utilizing their properties under the CRE through the conditioning argument; see details at ref. \cite{li2018asymptotic}.
Proposition \ref{prop:rem}(iii) involves careful analysis of the non-Gaussian distribution.  
% which itself is not directly related to the randomization-based asymptotics. 

\section{Permutational Berry-Esseen bounds}\label{sec:PBEB}

\subsection{Several univariate and multivariate permutational Berry-Esseen bounds}\label{sec:PBEB-theory}
Recently, permutational Berry-Esseen bounds (permutational BEBs, also called combinatorial BEBs) start to raise attention in randomization-based inference for experiments. Berry-Esseen bounds depict the distance between the sampling distribution of a statistic and a target, often normal, distribution. Theoretically speaking, it measures the convergence rate of central limit theorems. 
In general, the distance between two probability distributions is based on a class of metric of the following form:
% that is used to bound the 
\begin{align}
    d(\bbP_1, \bbP_2) = \sup_{h\in\cH} \lt| \int h \mathrm{d}\bbP_1 - \int h \mathrm{d}\bbP_2 \rt|.
\end{align}
In particular, Berry-Esseen bounds 
% target 
consider $\cH$ to be
the class of indicator functions over a family of sets. For univariate distributions, BEBs study the upper bound based on the Kolmogorov metric, where $\cH$ contains half-line indicator functions:
\begin{align}
    \sup_{t\in\bbR} \lt|\bbP_1\{X \le t\} - \bbP_2\{X \le t\}\rt|.
\end{align}
In the multivariate case, there are many choices of sets for different purposes, such as Euclidean balls \citep{bentkus2003dependence}, rectangular sets \citep{chernozhukov2017central}, measurable convex sets \citep{bentkus2005lyapunov, bhattacharya2010normal}, etc. 

Below we review some theoretical progresses of Permutational BEBs and their important application for analyzing randomized experiments in finite populations. 

\subsubsection{Univariate case} \label{sec:PBEB-theory-univariate}

We consider the univariate linear permutational statistics in \eqref{eqn:sips} and adopt the notation from Section \ref{sec:linear_perm_clt}. We will summarize BEB results for $\Gamma$ upon standardization. The standardized version of $\Gamma$ can be expressed as
\begin{align}
    {\Var{\Gamma}}^{-1/2}(\Gamma - \E{\Gamma}) = \frac{\sum_{i=1}^N \tilde{M}_N(i,\pi(i))}{\{\frac{1}{N-1}\sum_{i,j=1}^N \tilde{M}_N(i,j)^2\}^{1/2}} = \sum_{i=1}^N \check{M}_N(i, \pi(i)),
\end{align}
where
\begin{align}
    \check{M}(i,j) = \frac{\tilde{M}_N(i,j)}{\{(N-1)^{-1}\sum_{i,j=1}^N \tilde{M}_N(i,j)^2\}^{1/2}}. 
\end{align}
Therefore, without loss of generality, we assume the following condition. 
\begin{condition}[Normalizing $\Gamma$]\label{cond:normalizing-M}
 $\Gamma$ in Definition \eqref{eqn:sips} is defined with $M_N$ satisfying the following normalizing condition:
\begin{align}\label{eqn:suff-M}
    M_N(i, +) = M_N(+, j) = M_N(+, +) = 0, \text{ for all } i,j\in[N]; \quad \sum_{i,j\in[N]}M_N(i,j)^2 = N-1. 
\end{align} 
\end{condition}

\citet{von1976remainder} and \citet{ho1978lp} established some early results. \citet{bolthausen1984estimate} applied one version of Stein's method \citep{stein1972bound} to establish the following result requiring only conditions concerning the third moment of the matrix $M_N$. 
\begin{proposition}[Main theorem of \citet{bolthausen1984estimate}]\label{prop:beb-B}
Assume Condition \ref{cond:normalizing-M}. There exists some universal constant $C>0$, such that
\begin{align}\label{eqn:beb-B}
\sup_{t\in\bbR} |\bbP\{\Gamma \le t\} - \Phi(t)| \le 
 CN^{-1}\sum_{i,j\in[N]}|{M}_N(i,j)|^3.
\end{align}
\end{proposition}

The upper bound on the right-hand-side of 
\eqref{eqn:beb-B} achieves the rate of 
$O(N^{-1/2})$ if 
\begin{align}\label{eqn:beb-B-suff}
{N^{-1/2}\sum_{i,j\in[N]}|{M}_N(i,j)|^3} = O(1).
\end{align}
\citet{von1976remainder} imposed the following boundedness condition, which
is sufficient for \eqref{eqn:beb-B-suff}:  
% and thus \eqref{eqn:beb-B}:
\begin{align}
    \sup_{i,j\in[N]}|{M}_N(i,j)| = O(N^{-1/2}).
\end{align}
As a side note, the above boundedness condition is also sufficient for the BEB in ref. \cite{ho1978lp} to achieve the $O(N^{-1/2})$ convergence rate. \citet[Chapter 6.1]{chen2011normal} presented a thorough discussion about the univariate permutational BEB. Proposition \ref{prop:beb-B} is very helpful for analyzing the finite-sample quality of normal approximation for linear estimators. In Section \ref{sec:apply-PBEB}, we provide an example of using Proposition \ref{prop:beb-B} to analyze CREs with possibly varying group sizes and diverging treatment levels. 

% \textbf{Multivariate case.} 
\subsubsection{Multivariate case}\label{sec:PBEB-theory-multivariate}
% Furthermore, 
We now consider the multivariate linear permutational statistics \eqref{eqn:mv-sips} in Section \ref{sec:mv-sips}. For the ease of presentation, we focus on results upon standardization. 
Specifically, 
\citet[Lemma S2]{shi2022berry} proved that the standardized version of \eqref{eqn:mv-sips} can still be written as a multivariate linear permutational statistics with a different set of $\check{M}_{N,h}$'s:
\begin{align}
    \Var{\Gamma}^{-1/2}(\Gamma - \E{\Gamma}) = \lt(\sum_{i=1}^N \check{M}_{N,1}(i,\pi(i)), \dots, \sum_{i=1}^N \check{M}_{N,H}(i,\pi(i))\rt)^\top, 
\end{align}
where $\check{M}_{N,h}$'s satisfy the following normalizing conditions:
\begin{gather}
    \check{M}_{N,h}(i, +) = \check{M}_{N,h}(+, j) = \check{M}_{N,h}(+, +) = 0, \text{ for all } i,j\in[N] \text{ and } h\in[H]; \label{eqn:center-Mh} \\
    \sum_{i,j\in[N]} \check{M}_{N,h}(i,j)^2 = N-1, \text{ for all } h\in[H]; \label{eqn:unit-Mh}\\
    \sum_{i,j\in[N]} \check{M}_{N,h}(i,j) \check{M}_{N,h'}(i,j) = 0, \text{ for all } h\neq h'\in[H]. \label{eqn:orthogonal-Mh}
\end{gather}
$\check{M}_{N,h}$'s can be constructed from $M_{N,h}$'s by performing the centering step as in \eqref{eqn:tildeM} then applying a linear combination using  the matrix $ \Var{\Gamma}^{-1/2} $. Therefore, without loss of generality, we assume the following condition. 
\begin{condition}[Normalizing $\Gamma$ in the multivariate case]\label{cond:normalizing-Mh}
    $\Gamma$ in \eqref{eqn:mv-sips} is defined with $M_{N,h}$'s satisfying the normalizing conditions in \eqref{eqn:center-Mh}, \eqref{eqn:unit-Mh} and \eqref{eqn:orthogonal-Mh}, which guarantees that $\E{\Gamma} = 0_H$ and $\Var{\Gamma} = I_{H}$. 
\end{condition}

\citet{bolthausen1993rate} extended 
the univariate 
% case 
result in Proposition \ref{prop:beb-B}
to the multivariate, possibly nonlinear setting. 
In particular, for the multivariate linear case, \citet[Theorem 1]{bolthausen1993rate} established the following BEB. 
\begin{proposition}
    % {\lxr clarify $M$}
Let $\cA$ be the family of all measurable convex sets. Let $\Gamma_Z$ be a random Gaussian vector that follows $\cN(0_H, I_H)$. Assume Condition \ref{cond:normalizing-Mh}. Then there exists a constant $C_H$ that depends only on the dimension $H$ such that % \inlinemod{red}{What is $\Gamma_Z$?}
\begin{align}\label{eqn:beb-BG}
    \sup_{A\in\cA}|\bbP\{\Gamma \in A\} - \bbP\{\Gamma_Z \in A\}|\le   \frac{C_H}{N}\sum_{i,j\in[N]}\lt(\sum_{h=1}^H M_{N,h}(i,j)^2\rt)^{3/2}. 
\end{align}
\end{proposition}

The BEB in 
\eqref{eqn:beb-BG} covers the univariate case 
in Proposition \ref{prop:beb-B} as a special case with $H=1$.
% \eqref{eqn:beb-B} if $H = 1$. 
% The constant in \eqref{eqn:beb-BG}, $C_H$, depends on the dimension $H$. 
However, \citet{bolthausen1993rate} did not give a closed form expression for $C_H$, whose dependence on the dimension $H$ is unknown.  
\citet{raic2015multivariate} conjectured the following result:
\begin{align}\label{eqn:raic}
    \sup_{A\in\cA}|\bbP\{\Gamma\in A\} - \bbP\{\Gamma_Z\in A\}| \le C\frac{H^{1/4}}{N}\sum_{i,j\in[N]}\lt(\sum_{h\in[H]}M_{N,h}(i,j)^2\rt)^{3/2},
\end{align}
where $C_H$ can be an absolute constant that does not depend on the dimension $H$.
However, no formal proof is provided by the author. \citet{chatterjee2007multivariate} made one step forward to reveal the dimensional dependence using Stein's method with multivariate exchangeable pairs. In \citet[Section 3.2]{chatterjee2007multivariate}, the authors established a bound for the following distance: 
\begin{align}
    \sup_{g\in C^2(\bbR^H)} \lt|\E{g(\Gamma)} - \E{g(\Gamma_Z)}\rt|,
\end{align}
where $C^2(\bbR^H)$ stands for the class of $2$-times continuously differentiable functions on $\bbR^H$. We state below a special case of \citet{chatterjee2007multivariate}'s result. 
\begin{proposition}
    Under Condition \ref{cond:normalizing-Mh} and the condition of bounded entries:
\begin{align}\label{eqn:mv-beb-suff}
    \sup_{i,j\in[N], h\in[H]}|M_{N,h}(i,j)| = O(N^{-1/2}),
\end{align}
we have % {\lxr scale of $g$}
\begin{align}\label{eqn:beb-cm}
    \sup_{g\in C^2(\bbR^H)} \lt|\E{g(\Gamma)} - \E{g(\Gamma_Z)}\rt| = O\lt(\frac{H^3}{N^{1/2}}\rt).
\end{align}
\end{proposition}
Nevertheless, \eqref{eqn:beb-cm} does not translate directly into a Berry-Esseen bound under the Kolmogorov metric because the indicator functions are not members of $C^2(\bbR^H)$. \citet[Theorem S2]{shi2022berry} made use of one key result established by \citet{fang2015rates} regarding Stein's coupling and established the following multivariate permutational BEB with explicit dependence on the dimension. 
\begin{proposition}\label{prop:beb-S}
    Under Condition \ref{cond:normalizing-Mh} and the condition of bounded entries \eqref{eqn:mv-beb-suff}, we have 
\begin{align}
    \sup_{A\in\cA}|\bbP\{\Gamma \in A\} - \bbP\{\Gamma_Z \in A\}| = O\lt(\frac{H^{13/4}}{N^{1/2}}\rt).
\end{align}
\end{proposition}
Proposition \ref{prop:beb-S} is also useful for analyzing the finite-sample performance of many non-linear permutational statistics. \citet{shi2022berry} used Proposition \ref{prop:beb-S} to obtain a BEB for quadratic forms of a multi-dimensional estimator for causal effects in CRE, which builds up the ground for Wald-type inference. See Appendix E of \citet{shi2022berry} for more discussion. 

\subsection{Application of permutational BEBs to randomization-based inference}\label{sec:apply-PBEB}

In this section, we present several applications of permutational BEBs in randomization-based inference. 

% \textbf{CRE with possibly varying group sizes and diverging treatment levels.} 

\subsubsection{Completely randomized experiments with possibly varying group sizes and diverging treatment levels}

Many classical experiments only involve a small number of treatment levels. For example, classical factorial experiments typically include a small number of factors (like $K \le 5$) \citep{wu2011experiments}. However, many modern experiments involve a much larger number of treatment levels and units due to the need for analyzing more complex relations as well as the development of experimentation technologies. For example, in political science, powered by the development of computers and web-based technology, conjoint survey experiments \citep{caughey2019item, hainmueller2014causal, zhirkov2022estimating} (as a special type of factorial experiments) are very popular for analyzing the effects of many factors together and answering complex causal questions. In ref. \cite{zhirkov2022estimating}, the author investigated an experiment examining the impact of six different attributes of immigrants on public support for their admission to the United States. \citet{caughey2019item} studied the impact of twelve ($K = 12$) personal traits on citizens' preference for U.S. presidential candidates. A large number of treatment levels pose new challenges to the analysis of randomized experiments and call for new methodological and theoretical developments.  \citet{shi2022berry} and \citet{shi2023forward} discussed general CREs where the number of treatment levels $Q$ and the treatment group sizes $N_q$'s follow a variety of asymptotic regimes beyond the classical setup. Table \ref{tab:many-settings} presents several possible regimes that are of interest both technically and practically. 
\begin{table}[!ht]
\centering
\caption{Theoretical results for multi-level experiments under the randomization model, originally from Table 1 of ref. \cite{shi2022berry}. The column title ``$Q$'' and ``$N_q$'' stand for the number of treatment levels and the number of replications within the treatment levels, respectively. The last column summarizes how well each of the regimes is studied in the literature regarding CLT, variance estimation and BEB.}
\label{tab:many-settings}
\begin{threeparttable}
\begin{tabular}{M{1.0cm}M{1.0cm}M{3cm}M{7cm}}
 Regime & $Q$    & $N_q$                 & CLT, variance estimation, and BEB              \\ \hline
   (R1) &  Small &   Large & CLT and variance estimation established; no BEB  \\  
 (R2)       & Large & Large                    & Seems similar to (R1) but not studied \\  
 (R3)        & Large & Small but $N_q\ge$ 2 & Not studied    \\  
   (R4) &  Large &  $N_q = 1$ & Not studied; variance estimation is nontrivial \\ \hline
   (R5) &   \multicolumn{2}{c}{{Mixture of the above}} &  Not studied \\ 
\end{tabular}%
% \begin{tablenotes}
% \small
%        \item Note: The column title ``$Q$'' and ``$N_q$'' stand for the number of treatment levels and the number of replications within the treatment levels, respectively.
% \end{tablenotes}
\end{threeparttable}
\end{table}

Most of the regimes in Table \ref{tab:many-settings} are less visited by literature and lacks  scientific justification. \citet{shi2022berry} utilized permutational BEBs to characterize 
the normal approximation for sampling distributions of statistics
% the finite sample distance of the sampling statistics 
in general CREs, and managed to present a unified discussion of all the regimes listed in Table \ref{tab:many-settings}. We elaborate on the usage of permutational BEBs with a canonical example in factorial experiments from \citet{shi2022berry}. 

In a $2^K$ factorial design with $K$ binary factors, there are  $Q=2^K$ possible treatment levels. Index the potential outcomes $Y_i(q)$'s also as $Y_i(z_1, \ldots, z_K)$'s, where $q = 1,\ldots, Q$ and $z_1,\ldots, z_K = 0,1$. The parameter of interest $\tau = F^\top \overline{Y}(\cdot)  $ may consist of a subset of factorial effects, where $F$ is a contrast matrix with orthogonal columns and entries of $\pm (Q/2)^{-1}$;  %{\color{red} [R1: Example 1: ``$\pm Q^{-1}$ entries" $\to$ ``entries of $\pm Q^{-1}$".]} {\color{orange} entries of $\pm Q^{-1}$}; 
see \citet{dasgupta2015causal} for precise definitions of main effects and interactions. The factorial design is called \textit{nearly uniform} if the sizes of each arm, $N_q$'s, are approximately of the same order. More rigorously, we assume that
there exists a positive integer $N_0 > 0$ and absolute constants $\underline{c} \le \overline{c}$ such that $ N_q = c_q  {N}_0$ with  $\underline{c}\le c_q\le \overline{c}$, for all $q=1,\ldots, Q$. 
% Basically such a setup covers 
Such a setup can cover many cases in 
regimes (R1)-(R4) in Table \ref{tab:many-settings}. 
\citet{shi2022berry} established the following result for the plug-in estimator $\htau = F^\top \hY(\cdot)$:

\begin{proposition}[\citet{shi2022berry}, Example 6, nearly uniform factorial design]\label{prop:beb-factorial}
Consider a nearly uniform $2^K$ factorial experiment. Let $\ttau = \Var{\htau}^{-1/2}(\htau - \tau)$ be the standardized version of $\htau$. Let $F\in\bbR^{Q\times H}$ with $H = K+K(K-1)/2={K(K+1)}/{2}$ be the contrast matrix for all main effects and two-way interactions. Recall the definition of $S(q,q)$ from \eqref{eq:S_qq}. Under some mild regularity conditions, we have
%\begin{align} \label{eqn:factorial-BE}
% \sup_{b\in\bbR^H,\|b\|_2 =1}\sup_{t\in\bbR}\lt|\bbP\{b^\top \tilde{\gamma}  \le t\} - \Phi(t)\rt|  \le  C\sigma_F   \frac{  \max_{i\in[N],q\in[Q]}|Y_i(q)-\overline{Y}(q)|}{\{\min_{q\in[Q]} S(q,q)\}^{1/2} } \sqrt{\frac{H}{N}} .
%\end{align}
\begin{align} \label{eqn:factorial-BE}
 \sup_{b\in\bbR^H,\|b\|_2 =1}\sup_{t\in\bbR}\lt|\bbP\{b^\top\ttau   \le t\} - \Phi(t)\rt|  
 \le  
 C
 \sigma_F     
 \frac{ \max_{q\in[Q],i\in[N]} |Y_i(q) - \oY(q)| }{ \{  \min_{q\in[Q]} S(q,q) \}^{1/2} }  \sqrt{\frac{K^2}{N}}   .
\end{align}
where $C>0$ is an absolute constant and $\sigma_F>0$ is certain constant related to the matrix $F$.
%One can see that under the condition
%\begin{align*}
%    \frac{\underline{c}^{-1}  \max_{i\in[N],q\in[Q]}|Y_i(q)-\overline{Y}(q)|}{\{\overline{c}^{-1} \min_{q\in[Q]} S(q,q)\}^{1/2} } = o(N_0^{1/2} Q^{1/2} K^{-1}),
%\end{align*}
%the upper bound in \eqref{eqn:factorial-BE} converges to 0.  
\end{proposition}
Proposition \ref{prop:beb-factorial} is established based on the permutational BEB from \citet{bolthausen1984estimate} (presented in Proposition \ref{prop:beb-B} in Section \ref{sec:PBEB-theory}). In particular, one can formulate $b^\top\ttau$ as a linear permutational statistic with a carefully defined matrix $M_N$, and apply Proposition \ref{prop:beb-B} to obtain a raw BEB. After taking supreme over all unit-norm vector $b$, the BEB can be simplified to the presented form \eqref{eqn:factorial-BE}. More technical details are provided in Appendix A of ref. \cite{shi2022berry}. Also, the BEB in \eqref{eqn:factorial-BE} is uniform in the linear coefficient vector $b\in\bbR^H$ with $\|b\|_2 = 1$. This uniformity results in the additional dependence in $K^2$ (or the dimension $H$). Intuitively with higher dimension $H$ the uniform bound becomes larger. From Proposition \ref{prop:beb-factorial}, we can obtain a sufficient condition for the upper bound \eqref{eqn:factorial-BE} to converge to 0, which implies a CLT for any one-dimensional linear transformation of $\ttau$. 
%\todo{is $H$ or $K$ fixed?} 
In addition, Proposition \ref{prop:beb-factorial} requires mainly the total sample size $N$ to be large enough, and therefore allows either a fixed number of treatment levels $Q$ and diverging replications $N_0$, or a diverging $Q$ with limited replications $N_0$.
\citet{shi2022berry} also established design-based properties of Wald-type inference under general CREs, which utilizes multivarite permutational BEBs such as Proposition \ref{prop:beb-S}.

\subsubsection{Rerandomization with diminishing covariate imbalance and diverging number of covariates}

\citet{li2018asymptotic} studied the asymptotic theory of rerandomization with a fixed covariate imbalance threshold that does not vary with the sample size, as discussed in Sections \ref{sec:rerandomization} and \ref{sec:rem_fix_a}. 
The theory there suggests that the smaller the threshold, the more improvement we can gain from rerandomization over the complete randomization. 
Although intuitive, such a conclusion is not precise. 
When the covariate balance criterion is too stringent, there may be no acceptable assignments, 
and, even if there are acceptable ones, the asymptotic approximation may work poorly due to the small and even diminishing acceptance probability, i.e., the probability that a complete randomization is acceptable under rerandomization. 
Specifically and technically, 
the derivation for properties of rerandomization is through analyzing conditional distributions under the CRE, 
which will involve the acceptance probability in the denominator.
The resulting asymptotic analysis will then encounter a ratio between two quantities of order $o(1)$ when we allow the acceptance probability (or the imbalance threshold) to diminish with the sample size. 
In such cases, BEBs are crucial for conducting asymptotic analysis.

In the context of 
% complete randomization or equivalently 
simple random sampling, 
\citet{wang2022rerandomization} derived a multivariate BEB for the sample average using \citet{hajek1960limiting}'s coupling and the BEB for sums of independent random vectors \citep{raic2019multivariate} 
with explicit dependence on the dimension. 
The bound, although weaker than that implied by the conjecture in Raic \citep{raic2019multivariate}, 
is sufficient for studying rerandomization with diminishing covariate imbalance threshold (or equivalently acceptance probability) and diverging number of covariates. 
With the derived BEBs, \citet{wang2022rerandomization} presented the following asymptotic theory for ReM, which is stronger than Proposition \ref{prop:rem}. 
We adopt the same notation from Section \ref{sec:rem_fix_a}, 
and denote the covariance imbalance threshold by $a_n$ and the number of covariates by $K_n$, allowing them to vary with the sample size $N$. 
Let $r_1 = N_1/N$, $r_0 = N_0/N$, $u_i = (r_0 \cdot Y_i(1) + r_1 \cdot Y_i(0), X_i^\top)^\top$, 
$\bar{u}$ and $S_u^2$ be the finite population average and covariance of $u_i$'s,
and 
\begin{align*}
    \gamma_N & \equiv 
    \frac{(K_N+1)^{1/4}}{\sqrt{N r_1 r_0}}
    \frac{1}{N} \sum_{i=1}^N 
    \| S_u^{-1} (u_i - \bar{u}) \|_2^3, 
\end{align*}
where
% $K_n$ denotes the number of covariates used in rerandomization and 
$S_u^{-1}$ is the inverse of the positive semidefinite square root of $S_u^2$. 
We have the following BEB under ReM.
\begin{proposition}[\citet{wang2022rerandomization}, Theorems 1 and 3]\label{prop:ReM-BEB}
    % Assume that, 
    As $N\rightarrow \infty$, 
    if 
    $\gamma_N \rightarrow 0$ and $p_N/\gamma_N^{1/3} \rightarrow \infty$ with $p_N \equiv \Pr(\chi^2_{K_N}\le a_N)$, then 
    \begin{align*}
        \sup_{c\in \mathbb{R}}
        \big| 
        \Pr\{ \Var{\htau}^{-1/2}(\hat{\tau}-\tau) \le c \mid M\le a_N \}
        - 
        \Pr(
        \sqrt{1-R^2} \ \varepsilon_0 + \sqrt{R^2}\ L_{K_N, a_N} \le c
        )
        \big|
        \rightarrow 0, 
    \end{align*}
    where $\tau$ is the true average treatment effect.
\end{proposition}

\citet{wang2022rerandomization} further studied additional conditions such that the constrained Gaussian random variable $L_{K_N, a_N}$ becomes ignorable as $N\rightarrow \infty$, under which $\Var{\htau}^{-1/2}(\hat{\tau}-\tau)$ can asymptotically follow the Gaussian distribution $\mathcal{N}(0, 1-R^2)$ under rerandomization. 
This is the ideally optimal precision that one can expect under rerandomization, since the remaining variation is due to the part of potential outcomes that cannot be linearly explained by the covariates. 
Moreover, 
the Gaussian asymptotic distribution is the same as that of Lin's regression-adjusted estimator under the CRE. 
Intuitively, 
rerandomization and 
covariate adjustment are dual of each other, where the former is at the design stage while the latter is at the analysis stage.   
\citet{wang2022rerandomization} further proposed large-sample valid confidence intervals for the average treatment effect under rerandomization. 

% {\color{red} For Xinran to write}

\section{Extensions} 

\citet{neyman1923application} has motivated many important extensions for the design and analysis of randomized experiments, and the technical tools regarding permutations have been evolving during the past century. In this section, we discuss some other extensions beyond \citet{neyman1923application}.
% \todo{better way} 

% beyond the \citet{neyman1923application}'s proposal. 

\subsection{Other randomized experiments}

In this section we discuss several other widely used and studied randomized experiments, beyond \citet{neyman1923application}'s focus on the CRE. 

% hanzhong's paper in factorial experiment
\subsubsection{Stratified (block) randomized experiments}\label{sec:block}
Stratified randomized experiments (SRE) have been used widely in many fields, including agricultural study \citep{petersen1994agricultural}, biomedical study \citep{goldner1971effects}, social science \citep{chong2016iron}, etc. A SRE combines several different CREs according to the levels of a stratum indicator. Concretely speaking, consider an experiment with $K$ strata. Denote the number and proportion of units in stratum $k$ as $N_{[k]}$ and $\pi_{[k]} = N_{[k]}/N$, respectively,  where $k = 1,\dots, K$. Within stratum $k$, $N_{[k]1}$ units are randomized to receive treatment and the remaining $N_{[k]0} = N_{[k]} - N_{[k]1}$ units are assigned to control. Across strata, the randomization is conducted independently. The treatment assignment distribution is uniform over all possible randomizations. 

Analogous to CRE, in SRE, for unit $i$ in stratum $k$, we have potential outcomes $Y_{ki}(1)$ and $Y_{ki}(0)$ and individual causal effect $\tau_{ki} = Y_{ki}(1)-Y_{ki}(0)$. For stratum $k$, we have stratum-specific average causal effect
\begin{align}
    \tau_{[k]} = N_{[k]}^{-1} \sum_{i = 1}^{N_{[k]}} \tau_{ki}. 
\end{align}
The {\lxr overall} average causal effect is
\begin{align}
    \tau = N^{-1} \sum_{k=1}^K\sum_{i=1}^{N_{[k]}}\tau_{ki} = \sum_{k=1}^K \pi_{[k]}\tau_{[k]}. 
\end{align}
which is also a weighted average of the stratum-specific average causal effects. For Neyman-type analysis, a point estimator can be obtained by taking a weighted average of stratum-specific difference-in-means estimators:
\begin{align}\label{eqn:tau-S}
    \htau_\scS = \sum_{k=1}^K \pi_{[k]}\htau_{[k]},
\end{align}
where $\htau_{[k]}$ is the difference-in-means estimator for stratum $k$. It has variance
\begin{align*}
    \Var{\htau_\scS} = \sum_{k=1}^K \pi_{[k]}^2 \Var{\htau_{[k]}},
\end{align*}
which motivates the variance estimator 
\begin{align*}
    \hV_\scS = \sum_{k=1}^K \pi_{[k]}^2 \lt\{\frac{\hS_{[k]}^2(1)}{N_{[k]1}} + \frac{\hS_{[k]}^2(0)}{N_{[k]0}}\rt\}, 
\end{align*}
with $ \hS_{[k]}^2(1) $ and $ \hS_{[k]}^2(0) $ being the stratum-specific sample variances 
% estimators 
for the treatment and control arms.  
A Wald-type confidence interval can then be constructed for $\tau$.

% A confidence interval for $\tau$ is
% \begin{align}\label{eqn:CI-SRE}
%     \lt[\htau_\scS - z_{\alpha/2} \sqrt{\hv_\scS}, ~ \htau_\scS + z_{\alpha/2}\sqrt{\hv_\scS}\rt]. 
% \end{align}

Under certain regularity conditions, the point estimator \eqref{eqn:tau-S} is asymptotically normal and Wald-type inference is proved to be asymptotically valid (see, for example, refs \cite{Bickel1984, liu2020regression}). The random assignment mechanism requires studying 
% a mixture 
a convolution
of independent permutational distributions, which  motivates new theoretical tools. When the total number of strata $K$ is small and the sizes of the strata are large, the permutational CLTs 
%(e.g., Hoeffding's condition \eqref{eqn:sips-hoeffding})
play a central role in the analysis. When $K$ is large and the sizes of the strata are small, CLTs for independent summations play a crucial role instead. With a mixture of large and small strata, there are also theoretical results in the literature; see, for example, refs \cite{Bickel1984, liu2020regression}. Moreover, \citet{liu2022randomization} and \citet{wang2023rerandomization} further investigated covariate adjustment and rerandomization in SREs. 

\subsubsection{Matched-pairs experiments}
The matched-pairs experiment (MPE) is another popular experimental design in practice \citep{fisher1935design, ball1973reading, imai2008variance}. The MPE is the most extreme version of the SRE with only one treated unit and one control unit within each stratum, which is called a \textit{pair}. We can adopt the notations for the SRE to define potential outcomes, causal effects, stratum-specific difference-in-means estimator (denoted again as $\htau_{[k]}$) and the aggregated difference-in-means estimator (denoted as $\htau_\scM$) in the MPE. However, 
the variance estimation strategy discussed in Section \ref{sec:block} is no longer applicable for the MPE, since it implicitly requires at least two treated and control units within each matched set so that we can calculate the stratum-specific sample variances.
\citet{imai2008variance} proposed the following variance estimator by instead considering the sample variance of the stratum-specific difference-in-means estimators:
% variance estimation is nontrivial in MPE because only one unit receives treatment/control within each stratum.\todo{I will edit this} Many works in the literature discussed variance estimation strategies in such scenarios; see for example \cite{imai2008variance, fogarty2018mitigating, pashley2021insights}. Here we present a variance estimator proposed by \citet{imai2008variance}:
\begin{align*}
    \hV_\scM = \frac{1}{n(n-1)}\sum_{k=1}^n (\htau_{[k]} - \htau_\scM)^2, 
\end{align*}
and he showed that it is conservative in expectation for the true variance of $\htau_\scM$. 
% \citet{imai2008variance} showed its conservativeness for estimating $\Var{\htau_\scM}$. 
We can then construct the Wald-type confidence interval
\begin{align}
    \lt[\htau_\scM - z_{\alpha/2} {\hV_\scM}^{1/2}, ~ \htau_\scM + z_{\alpha/2} {\hV_\scM}^{1/2}\rt],
\end{align}
which can be 
% proves to be 
asymptotically valid under certain regularity conditions.
% for Type I error control. 
Moreover, regression adjustment can be applied to improve efficiency when baseline covariates are available, as shown in ref. \cite{fogarty2018regression}.

In general stratified experiments with possibly one treated or one control unit in some strata, 
\citet{fogarty2018mitigating} and \citet{pashley2021insights} discussed general strategies to conservatively estimate the variance of the aggregated difference-in-means estimator.

\subsubsection{Cluster randomized experiments}
Cluster randomized experiments are widely used due to their logistical convenience and policy relevance. In a cluster randomized experiment, the treatment is assigned at the cluster level instead of the individual level. Consider a study with $N$ units and $M$ clusters. Cluster $i$ has
$n_i$ units $(i = 1, \dots, M)$. Let $(i, j)$ index the $j$-th
unit within cluster $i$ for $i = 1, \dots , M$ and $j = 1, \dots, n_i$. 
% Unit $(i, j)$ has covariates $x_{ij}$, and cluster $i$ has covariates $c_i$. 
The experimenter randomly assigns 
% $eM$ 
$M_1$
clusters to receive the treatment and 
% $(1-e)M$ 
$M_0$
clusters to receive the control, where 
% $0 < e < 1$ is a fixed number denoting the proportion of treated
% clusters. 
{\lxr $M_1$ and $M_0$ are fixed positive integers satisfying $M_1+M_0=M$.}
Let $Z_i$ be the treatment indicator for cluster $i$ and $Z_{ij}$ be the treatment indicator for unit
$(i, j)$. In a cluster-randomized experiment, units within a cluster receive identical treatment levels.
So if cluster $i$ receives treatment, then $Z_{ij} = Z_i = 1$ for all $j$. If cluster $i$ receives control, then $Z_{ij} = Z_i = 0$ for all $j$. 
Let $Y_{ij}(1)$ and $Y_{ij}(0)$ be the potential outcomes under treatment
and control, respectively, for unit $(i, j)$. The observed outcome is 
% a function of the treatment indicator and potential outcomes: 
then $Y_{ij} = Z_{ij}Y_{ij}(1)+(1 - Z_{ij})Y_{ij}(0)$.
The average treatment effect over all units is
\begin{align}
    \tau = \frac{1}{N}\sum_{i=1}^M\sum_{j=1}^{n_i} \{Y_{ij}(1) - Y_{ij}(0)\}. 
\end{align}
There are different strategies for inferring $\tau$, including individual-level estimators
%\todo{different weight?} 
and cluster-level estimators \citep{su2021model}, both enjoying desirable asymptotic properties implied by permutational CLTs. We refer 
interested readers to a collection of works on analyzing cluster-randomized experiments \citep{abadie2023should, middleton2015unbiased, su2021model, lu2023design, schochet2022design, athey2017econometrics}.

% \subsection{Other causal estimands}
% quantile treatment effects 

\subsection{Some technical aspects for permutations}\label{sec:technical-aspects}

Permutation is a core element in the design and analysis of randomized experiments, 
and its development, such as CLTs and BEBs, has involved many technical tools including moment matching \citep{hoeffding1951combinatorial}, coupling \citep{hajek1961some, hajek1968asymptotic}, Stein's method \cite{stein1972bound, bolthausen1984estimate, chatterjee2007multivariate, chen2011normal}, etc. 
In this section we discuss some technical aspects for analyzing permutation-related problems.

\subsubsection{Hajek's coupling}
Hajek's coupling is one technique developed in ref. \cite{hajek1960limiting} for proving central limit theorems in simple random sampling. The idea is based on constructing a coupling between simple random sampling and Bernoulli random sampling so that CLTs for i.i.d. sampling can be applied. \citet{wang2022rerandomization} utilized Hajek's coupling together with a  multivariate BEB for sum of independent random vectors \citep{raic2019multivariate} to study rerandomization with diminishing covariate imbalance. The techniques are useful to establish theories for a wide range of permutational statistics; see, e.g., refs \cite{hajek1968asymptotic, liu2022randomization}. 
% questions.

\subsubsection{Double and multiple permutations}
Nowadays, there are many new variants of permutations in designing randomized experiments. For example,  \citet{fredrickson2019permutation} and \citet{chen2017new} discussed permutation and randomization tests for analyzing network data. \citet{d2016causal} and  \citet{deng2024unbiased} studied randomized experiments with dyadic outcomes, i.e., outcomes that measure the relationship between pairs of units. When randomization occurs at the unit level, the dyadic outcomes are in turn randomized with double permutations. \textit{Doubly indexed permutation statistics} (DIPS) is useful in these settings because the dyadic potential outcomes are functions of pairs of treatments, and the statistics for studying causal effects in these problems are generally represented as DIPS. \citet{bajari2021multiple, bajari2023experimental} proposed multiple randomization designs for marketplaces in which multiple populations interact and causal questions regarding interference are of particular interest. In terms of technical tools that are potentially useful for analyzing double or multiple permutations, \citet{chen2011normal, zhao1997error, reinert2007multivariate}, among others, used Stein's method \citep{stein1972bound} to study the asymptotic properties of DIPS. 

\subsubsection{Concentration inequalities}
Another technical tool that has been recognized by the literature is permutational/combinatorial concentration inequalities. \citet{bloniarz2016lasso} and \citet{lei2020regression} used permutational concentration inequalities to analyze regression adjustment in CREs when the dimension of the covariates is diverging. 
% \citet{bloniarz2016lasso} used permutational concentration inequalities to analyze regression adjustment in CREs with high dimensional covariates.  
It will be interesting to explore related potential research questions that involve delicate analysis of finite sample properties of permutational statistics and inspire the use of concentration inequalities. 
 
\section{Conclusion}

In this review, we revisited the fundamental contributions of
\citet{neyman1923application}'s seminal work regarding the introduction of potential outcomes, the promotion of physical randomization and the emphasis of repeated sampling properties of statistics over the randomization distribution. These contributions lay down the foundation for the design and analysis of randomized experiments. We also reviewed permutational central limit theorems and Berry--Esseen bounds in great detail,  
% (central limit theorems or CLTs and Berry--Esseen bounds or BEBs)  
and listed applications of these technical results in randomization-based inference.

Beyond what we have covered in the review, many research 
% questions are 
work are 
closely related to \citet{neyman1923application}. From a technical point of view, many theoretical tools are not fully covered in the discussion. For example, when analyzing stratified randomized experiments, we need central limit theorems and Berry-Esseen bounds that combine the independent permutational distributions \citep{liu2020regression, liu2022randomization}.
This is also closely related to Rosenbaum's sensitivity analysis for matched observational studies with biased permutations in each matched sets \citep{rosenbaum2002observational, Rosenbaum2000sep, wu2023sensitivity}.
As another example, for the design and analysis of adaptive experiments, a general martingale structure typically exists which requires a martingale central limit theorem or Berry-Esseen result \citep{hu2006theory, hall2014martingale}.  

From a practical point of view, many real-world examples can motivate the study of new designs, outcomes, assumptions, and causal estimands under the finite population framework. For example, interference among units is a common phenomenon in many experimental and observational studies. The study of interference {\lxr and peer influence} has motivated a lot of new designs and methods, such as designing and analyzing bipartite experiments \citep{harshaw2023design, zigler2021bipartite}, multiple randomization \citep{bajari2021multiple}, randomized experiments with network interference \citep{leung2022causal}, group formulation design \citep{li2019peer, basse2019randomization}, etc.
Another example is randomization with missing observations or covariates. \citet{zhao2022adjust} discussed several strategies for randomization-based inference with missing covariates, and \citet{zhao2024covariate} further studied covariate adjustment in randomized experiments with both missing outcomes and covariates. Censored survival outcomes are another type of missingness that occurs frequently in clinical trials. In these settings, to test the null hypothesis of no treatment effect for any unit, \citet{rosenbaum2002observational} proposed randomization tests for censored outcomes using a partial ordering, and \citet{zhang2005asymptotic} established asymptotic normality of the randomization distribution of the log-rank statistic. Both approaches require the assumption of identical potential censoring times under treatment and control. Recently \citet{li2023randomization} relaxed this assumption and proved that, under a Bernoulli randomized experiment, with non-informative i.i.d. censoring, the log-rank test is asymptotically valid for testing Fisher’s null hypothesis of no treatment effect on any unit.

At the same time,
there have been extensive progresses for analyzing treatment effects 
from a super-population perspective, and many of them share similar spirit as the randomization-based inference \citep{DingLiMiratrix2017}. 
% Below we name a few. 
% , there have been many progresses for better analyzing 
For example, \citet{Tsiatis2001} have suggested linear covariate adjustment with treatment-covariate interaction under a semiparametric model; 
see also refs \cite{Rv2010, Robins2005}.  
% for more general covariate adjustment that can be robust to model misspecification. 
% for model-agnostic covariate adjustment for treatment effect estimation. 
In the presence of censoring, 
there have been many works studying semiparametric estimation of treatment effect \citep{laan2003unified, rubin2007doubly, van2011targeted}, as well as covariate adjustment to improve inference efficiency \citep{hernandez2006randomized, lu2008improving, moore2009increasing}. 

% \citet{laan2003unified, moore2009increasing, rubin2007doubly, van2011targeted}, among others, discussed semiparametric estimation of  survival causal parameter. Many works also discussed covariate adjustment with time-to-event outcomes using Cox proportional hazards models; see \citet{hernandez2006randomized, lu2008improving}, etc.

\bigskip 
\bigskip 
\bigskip

\noindent\textbf{Acknowledgements:} We thank the reviewers for carefully reading the manuscripts and providing many constructive suggestions for improving the paper.

\bigskip 
\bigskip 
\bigskip

\noindent\textbf{Funding information:} X. L. is partly supported by the National Science Foundation under grant DMS-2400961.

\bigskip 
\bigskip 
\bigskip

\noindent\textbf{Author contribution:}
All authors have accepted responsibility for the entire content of this manuscript and consented to its submission to the journal, reviewed all the results, and approved the final version of the manuscript. LS and XL worked together to frame the structure of the review, collect related literature, and organize the presentation of the methodology and theory for the design and analysis of randomized experiments.

\bibliographystyle{abbrvnat}
\bibliography{ref}

\begin{thebibliography}{163}
\providecommand{\natexlab}[1]{#1}
\providecommand{\url}[1]{\texttt{#1}}
\expandafter\ifx\csname urlstyle\endcsname\relax
  \providecommand{\doi}[1]{doi: #1}\else
  \providecommand{\doi}{doi: \begingroup \urlstyle{rm}\Url}\fi

\bibitem[Abadie et~al.(2023)Abadie, Athey, Imbens, and Wooldridge]{abadie2023should}
A.~Abadie, S.~Athey, G.~W. Imbens, and J.~M. Wooldridge.
\newblock {When should you adjust standard errors for clustering?}
\newblock \emph{The Quarterly Journal of Economics}, 138\penalty0 (1):\penalty0 1--35, 2023.

\bibitem[Angrist et~al.(1996)Angrist, Imbens, and Rubin]{angrist1996identification}
J.~D. Angrist, G.~W. Imbens, and D.~B. Rubin.
\newblock Identification of causal effects using instrumental variables.
\newblock \emph{Journal of the American statistical Association}, 91\penalty0 (434):\penalty0 444--455, 1996.

\bibitem[Athey and Imbens(2017{\natexlab{a}})]{ATHEY201773}
S.~Athey and G.~W. Imbens.
\newblock {The Econometrics of Randomized Experiments}.
\newblock In A.~Banerjee and E.~Duflo, editors, \emph{{Handbook of Economic Field Experiments}}, volume~1, chapter~3, pages 73--140. North-Holland, Amsterdam, 2017{\natexlab{a}}.

\bibitem[Athey and Imbens(2017{\natexlab{b}})]{athey2017econometrics}
S.~Athey and G.~W. Imbens.
\newblock {The econometrics of randomized experiments}.
\newblock In \emph{Handbook of economic field experiments}, volume~1, pages 73--140. Elsevier, 2017{\natexlab{b}}.

\bibitem[Bajari et~al.(2021)Bajari, Burdick, Imbens, Masoero, McQueen, Richardson, and Rosen]{bajari2021multiple}
P.~Bajari, B.~Burdick, G.~W. Imbens, L.~Masoero, J.~McQueen, T.~Richardson, and I.~M. Rosen.
\newblock Multiple randomization designs, 2021.

\bibitem[Bajari et~al.(2023)Bajari, Burdick, Imbens, Masoero, McQueen, Richardson, and Rosen]{bajari2023experimental}
P.~Bajari, B.~Burdick, G.~W. Imbens, L.~Masoero, J.~McQueen, T.~S. Richardson, and I.~M. Rosen.
\newblock {Experimental design in marketplaces}.
\newblock \emph{Statistical Science}, 1\penalty0 (1):\penalty0 1--19, 2023.

\bibitem[Ball et~al.(1973)]{ball1973reading}
S.~Ball et~al.
\newblock {Reading with Television: An Evaluation of The Electric Company. A Report to the Children's Television Workshop. Volumes 1 and 2.}
\newblock 1973.

\bibitem[Bang and Robins(2005)]{Robins2005}
H.~Bang and J.~M. Robins.
\newblock Doubly robust estimation in missing data and causal inference models.
\newblock \emph{Biometrics}, 61:\penalty0 962--973, 2005.

\bibitem[Basse et~al.(2019)Basse, Ding, Feller, and Toulis]{basse2019randomization}
G.~Basse, P.~Ding, A.~Feller, and P.~Toulis.
\newblock Randomization tests for peer effects in group formation experiments.
\newblock \emph{arXiv preprint arXiv:1904.02308}, 2019.

\bibitem[Bauer et~al.(2008)Bauer, Sterba, and Hallfors]{bauer2008evaluating}
D.~J. Bauer, S.~K. Sterba, and D.~D. Hallfors.
\newblock {Evaluating group-based interventions when control participants are ungrouped}.
\newblock \emph{Multivariate behavioral research}, 43\penalty0 (2):\penalty0 210--236, 2008.

\bibitem[Bentkus(2003)]{bentkus2003dependence}
V.~Bentkus.
\newblock {On the dependence of the Berry--Esseen bound on dimension}.
\newblock \emph{Journal of Statistical Planning and Inference}, 113\penalty0 (2):\penalty0 385--402, 2003.

\bibitem[Bentkus(2005)]{bentkus2005lyapunov}
V.~Bentkus.
\newblock {A Lyapunov-type bound in $\bbR^d$}.
\newblock \emph{Theory of Probability \& Its Applications}, 49\penalty0 (2):\penalty0 311--323, 2005.

\bibitem[Bhattacharya and Rao(2010)]{bhattacharya2010normal}
R.~N. Bhattacharya and R.~R. Rao.
\newblock \emph{{Normal Approximation and Asymptotic Expansions}}.
\newblock SIAM, 2010.

\bibitem[Bickel and Freedman(1984)]{Bickel1984}
P.~J. Bickel and D.~A. Freedman.
\newblock {Asymptotic Normality and the Bootstrap in Stratified Sampling}.
\newblock \emph{The Annals of Statistics}, 12:\penalty0 470 -- 482, 1984.

\bibitem[Bloniarz et~al.(2016)Bloniarz, Liu, Zhang, Sekhon, and Yu]{bloniarz2016lasso}
A.~Bloniarz, H.~Liu, C.-H. Zhang, J.~S. Sekhon, and B.~Yu.
\newblock {Lasso adjustments of treatment effect estimates in randomized experiments}.
\newblock \emph{Proceedings of the National Academy of Sciences}, 113\penalty0 (27):\penalty0 7383--7390, 2016.

\bibitem[Bolthausen(1984)]{bolthausen1984estimate}
E.~Bolthausen.
\newblock {An estimate of the remainder in a combinatorial central limit theorem}.
\newblock \emph{Zeitschrift f{\"u}r Wahrscheinlichkeitstheorie und verwandte Gebiete}, 66\penalty0 (3):\penalty0 379--386, 1984.

\bibitem[Bolthausen and Gotze(1993)]{bolthausen1993rate}
E.~Bolthausen and F.~Gotze.
\newblock {The rate of convergence for multivariate sampling statistics}.
\newblock \emph{The Annals of Statistics}, pages 1692--1710, 1993.

\bibitem[Branson and Dasgupta(2020)]{BD20}
Z.~Branson and T.~Dasgupta.
\newblock Sampling-based randomised designs for causal inference under the potential outcomes framework.
\newblock \emph{International Statistical Review}, 88:\penalty0 101--121, 2020.

\bibitem[Branson et~al.(2016)Branson, Dasgupta, and Rubin]{branson2016improving}
Z.~Branson, T.~Dasgupta, and D.~B. Rubin.
\newblock Improving covariate balance in 2 k factorial designs via rerandomization with an application to a new york city department of education high school study.
\newblock \emph{The Annals of Applied Statistics}, pages 1958--1976, 2016.

\bibitem[Branson et~al.(2023)Branson, Li, and Ding]{Zach2023}
Z.~Branson, X.~Li, and P.~Ding.
\newblock {Power and sample size calculations for rerandomization}.
\newblock \emph{Biometrika}, 111:\penalty0 355--363, 2023.

\bibitem[Bruhn and McKenzie(2009{\natexlab{a}})]{Bruhn:2009}
M.~Bruhn and D.~McKenzie.
\newblock In pursuit of balance: Randomization in practice in development field experiments.
\newblock \emph{American Economic Journal: Applied Economics}, 1:\penalty0 200--232, 2009{\natexlab{a}}.

\bibitem[Bruhn and McKenzie(2009{\natexlab{b}})]{bruhn2009pursuit}
M.~Bruhn and D.~McKenzie.
\newblock In pursuit of balance: Randomization in practice in development field experiments.
\newblock \emph{American economic journal: applied economics}, 1\penalty0 (4):\penalty0 200--232, 2009{\natexlab{b}}.

\bibitem[Caughey et~al.(2019)Caughey, Katsumata, and Yamamoto]{caughey2019item}
D.~Caughey, H.~Katsumata, and T.~Yamamoto.
\newblock Item response theory for conjoint survey experiments.
\newblock Technical report, Working Paper, 2019.

\bibitem[Caughey et~al.(2023)Caughey, Dafoe, Li, and Miratrix]{CDLM21quantile}
D.~Caughey, A.~Dafoe, X.~Li, and L.~Miratrix.
\newblock Randomization inference beyond the sharp null: Bounded null hypotheses and quantiles of individual treatment effects.
\newblock \emph{Journal of the Royal Statistical Society, Series B (Statistical Methodology)}, in press, 2023.

\bibitem[Chatterjee and Meckes(2007)]{chatterjee2007multivariate}
S.~Chatterjee and E.~Meckes.
\newblock {Multivariate normal approximation using exchangeable pairs}.
\newblock \emph{arXiv preprint math/0701464v1}, 2007.

\bibitem[Chen and Friedman(2017)]{chen2017new}
H.~Chen and J.~H. Friedman.
\newblock {A new graph-based two-sample test for multivariate and object data}.
\newblock \emph{Journal of the American statistical association}, 112\penalty0 (517):\penalty0 397--409, 2017.

\bibitem[Chen et~al.(2011)Chen, Goldstein, and Shao]{chen2011normal}
L.~H. Chen, L.~Goldstein, and Q.-M. Shao.
\newblock \emph{{Normal Approximation by Stein's Method}}, volume~2.
\newblock Springer, 2011.

\bibitem[Chen et~al.(2023)Chen, Li, and Zhang]{chen2023role}
Z.~Chen, X.~Li, and B.~Zhang.
\newblock The role of randomization inference in unraveling individual treatment effects in clinical trials: Application to hiv vaccine trials.
\newblock \emph{arXiv preprint arXiv:2310.14399}, 2023.

\bibitem[Chernozhukov et~al.(2017)Chernozhukov, Chetverikov, and Kato]{chernozhukov2017central}
V.~Chernozhukov, D.~Chetverikov, and K.~Kato.
\newblock {Central limit theorems and bootstrap in high dimensions}.
\newblock \emph{The Annals of Probability}, 45\penalty0 (4):\penalty0 2309, 2017.

\bibitem[Chong et~al.(2016)Chong, Cohen, Field, Nakasone, and Torero]{chong2016iron}
A.~Chong, I.~Cohen, E.~Field, E.~Nakasone, and M.~Torero.
\newblock {Iron deficiency and schooling attainment in Peru}.
\newblock \emph{American Economic Journal: Applied Economics}, 8\penalty0 (4):\penalty0 222--255, 2016.

\bibitem[Cochran(1977)]{cochran1977sampling}
W.~G. Cochran.
\newblock \emph{{Sampling techniques}}.
\newblock John Wiley \& Sons, 1977.

\bibitem[Cohen and Fogarty(2020)]{cohen2020no}
P.~L. Cohen and C.~B. Fogarty.
\newblock No-harm calibration for generalized oaxaca-blinder estimators.
\newblock \emph{arXiv preprint arXiv:2012.09246}, 2020.

\bibitem[Cohen and Fogarty(2022)]{cohen2022gaussian}
P.~L. Cohen and C.~B. Fogarty.
\newblock Gaussian prepivoting for finite population causal inference.
\newblock \emph{Journal of the Royal Statistical Society Series B: Statistical Methodology}, 84\penalty0 (2):\penalty0 295--320, 2022.

\bibitem[Copas(1973)]{copas1973randomization}
J.~B. Copas.
\newblock {Randomization models for the matched and unmatched 2$\times$2 tables}.
\newblock \emph{Biometrika}, 60\penalty0 (3):\penalty0 467--476, 1973.

\bibitem[Cox(2009)]{cox2009randomization}
D.~Cox.
\newblock Randomization in the design of experiments.
\newblock \emph{International Statistical Review}, 77\penalty0 (3):\penalty0 415--429, 2009.

\bibitem[Dasgupta et~al.(2015)Dasgupta, Pillai, and Rubin]{dasgupta2015causal}
T.~Dasgupta, N.~S. Pillai, and D.~B. Rubin.
\newblock {Causal inference from $2^K$ factorial designs by using potential outcomes}.
\newblock \emph{Journal of the Royal Statistical Society: Series B}, 77:\penalty0 727--753, 2015.

\bibitem[David(1938)]{david1938limiting}
F.~David.
\newblock {Limiting distributions connected with certain methods of sampling human populations}.
\newblock \emph{Stat. Res. Mem}, 2:\penalty0 69--90, 1938.

\bibitem[Deng et~al.(2024)Deng, Li, Zhang, Wang, and Chen]{deng2024unbiased}
L.~Deng, Y.~Li, J.~Zhang, Y.~Wang, and C.~Chen.
\newblock Unbiased estimation for total treatment effect under interference using aggregated dyadic data.
\newblock \emph{arXiv preprint arXiv:2402.12653}, 2024.

\bibitem[DiCiccio and Romano(2017)]{diciccio2017robust}
C.~J. DiCiccio and J.~P. Romano.
\newblock {Robust permutation tests for correlation and regression coefficients}.
\newblock \emph{Journal of the American Statistical Association}, 112\penalty0 (519):\penalty0 1211--1220, 2017.

\bibitem[Ding(2017)]{ding2017paradox}
P.~Ding.
\newblock A paradox from randomization-based causal inference.
\newblock \emph{Statistical science}, pages 331--345, 2017.

\bibitem[Ding(2023)]{ding2023first}
P.~Ding.
\newblock {A First Course in Causal Inference}.
\newblock \emph{arXiv preprint arXiv:2305.18793}, 2023.

\bibitem[Ding and Dasgupta(2018)]{ding2018randomization}
P.~Ding and T.~Dasgupta.
\newblock A randomization-based perspective on analysis of variance: a test statistic robust to treatment effect heterogeneity.
\newblock \emph{Biometrika}, 105\penalty0 (1):\penalty0 45--56, 2018.

\bibitem[Ding et~al.(2017)Ding, Li, and Miratrix]{DingLiMiratrix2017}
P.~Ding, X.~Li, and L.~W. Miratrix.
\newblock Bridging finite and super population causal inference.
\newblock \emph{Journal of Causal Inference}, 5:\penalty0 20160027, 2017.

\bibitem[Ding et~al.(2019)Ding, Feller, and Miratrix]{DFM2019}
P.~Ding, A.~Feller, and L.~Miratrix.
\newblock Decomposing treatment effect variation.
\newblock \emph{Journal of the American Statistical Association}, 114:\penalty0 304--317, 2019.

\bibitem[D’Amour and Airoldi(2016)]{d2016causal}
A.~D’Amour and E.~Airoldi.
\newblock Causal inference for dyadic outcomes in social network analysis.
\newblock 2016.

\bibitem[Fang and R{\"o}llin(2015)]{fang2015rates}
X.~Fang and A.~R{\"o}llin.
\newblock Rates of convergence for multivariate normal approximation with applications to dense graphs and doubly indexed permutation statistics.
\newblock \emph{Bernoulli}, pages 2157--2189, 2015.

\bibitem[Fienberg and Tanur(1996)]{fienberg1996reconsidering}
S.~E. Fienberg and J.~M. Tanur.
\newblock {Reconsidering the fundamental contributions of Fisher and Neyman on experimentation and sampling}.
\newblock \emph{International Statistical Review/Revue Internationale de Statistique}, pages 237--253, 1996.

\bibitem[Fisher(1925)]{fisher1925statistical}
R.~A. Fisher.
\newblock \emph{{Statistical Methods for Research Workers}}.
\newblock Edinburgh by Oliver and Boyd, 1st edition, 1925.

\bibitem[Fisher(1935)]{fisher1935design}
R.~A. Fisher.
\newblock \emph{{The Design of Experiments}}.
\newblock Edinburgh, London: Oliver and Boyd, 1st edition, 1935.

\bibitem[Fisher and Mackenzie(1923)]{fisher1923studies}
R.~A. Fisher and W.~A. Mackenzie.
\newblock {Studies in crop variation. II. The manurial response of different potato varieties}.
\newblock \emph{The Journal of Agricultural Science}, 13\penalty0 (3):\penalty0 311--320, 1923.

\bibitem[Fogarty(2018{\natexlab{a}})]{fogarty2018mitigating}
C.~B. Fogarty.
\newblock On mitigating the analytical limitations of finely stratified experiments.
\newblock \emph{Journal of the Royal Statistical Society Series B: Statistical Methodology}, 80\penalty0 (5):\penalty0 1035--1056, 2018{\natexlab{a}}.

\bibitem[Fogarty(2018{\natexlab{b}})]{fogarty2018regression}
C.~B. Fogarty.
\newblock {Regression assisted inference for the average treatment effect in paired experiments}.
\newblock \emph{Biometrika}, 105:\penalty0 994--1000, 2018{\natexlab{b}}.

\bibitem[Fraser(1956)]{fraser1956vector}
D.~Fraser.
\newblock {A vector form of the Wald-Wolfowitz-Hoeffding theorem}.
\newblock \emph{The Annals of Mathematical Statistics}, pages 540--543, 1956.

\bibitem[Fredrickson and Chen(2019)]{fredrickson2019permutation}
M.~M. Fredrickson and Y.~Chen.
\newblock {Permutation and randomization tests for network analysis}.
\newblock \emph{Social Networks}, 59:\penalty0 171--183, 2019.

\bibitem[Freedman(2008{\natexlab{a}})]{freedman2008Aregression}
D.~A. Freedman.
\newblock {On regression adjustments to experimental data}.
\newblock \emph{Advances in Applied Mathematics}, 40\penalty0 (2):\penalty0 180--193, 2008{\natexlab{a}}.

\bibitem[Freedman(2008{\natexlab{b}})]{freedman2008Bregression}
D.~A. Freedman.
\newblock {On regression adjustments in experiments with several treatments}.
\newblock \emph{Annals of Applied Statistics}, 2:\penalty0 176--96, 2008{\natexlab{b}}.

\bibitem[Gastwirth et~al.(2000)Gastwirth, Krieger, and Rosenbaum]{Rosenbaum2000sep}
J.~L. Gastwirth, A.~M. Krieger, and P.~R. Rosenbaum.
\newblock Asymptotic separability in sensitivity analysis.
\newblock \emph{Journal of the Royal Statistical Society: Series B}, 62:\penalty0 545--555, 2000.

\bibitem[Goldner et~al.(1971)Goldner, Knatterud, and Prout]{goldner1971effects}
M.~G. Goldner, G.~L. Knatterud, and T.~E. Prout.
\newblock Effects of hypoglycemic agents on vascular complications in patients with adult-onset diabetes: Iii. clinical implications of ugdp results.
\newblock \emph{JAMA}, 218\penalty0 (9):\penalty0 1400--1410, 1971.

\bibitem[Guo and Basse(2021)]{guo2021generalized}
K.~Guo and G.~Basse.
\newblock {The generalized Oaxaca-Blinder estimator}.
\newblock \emph{Journal of the American Statistical Association}, pages 1--13, 2021.

\bibitem[Hainmueller and Hopkins(2015)]{hainmueller2015hidden}
J.~Hainmueller and D.~J. Hopkins.
\newblock {The hidden American immigration consensus: A conjoint analysis of attitudes toward immigrants}.
\newblock \emph{American journal of political science}, 59\penalty0 (3):\penalty0 529--548, 2015.

\bibitem[Hainmueller et~al.(2014)Hainmueller, Hopkins, and Yamamoto]{hainmueller2014causal}
J.~Hainmueller, D.~J. Hopkins, and T.~Yamamoto.
\newblock {Causal inference in conjoint analysis: Understanding multidimensional choices via stated preference experiments}.
\newblock \emph{Political analysis}, 22\penalty0 (1):\penalty0 1--30, 2014.

\bibitem[H{\'a}jek(1960)]{hajek1960limiting}
J.~H{\'a}jek.
\newblock {Limiting distributions in simple random sampling from a finite population }.
\newblock \emph{Publications of the Mathematical Institute of the Hungarian Academy of Sciences}, 5:\penalty0 361--374, 1960.

\bibitem[H{\'a}jek(1961)]{hajek1961some}
J.~H{\'a}jek.
\newblock {Some extensions of the Wald-Wolfowitz-Noether theorem}.
\newblock \emph{The Annals of Mathematical Statistics}, pages 506--523, 1961.

\bibitem[H{\'a}jek(1968)]{hajek1968asymptotic}
J.~H{\'a}jek.
\newblock {Asymptotic normality of simple linear rank statistics under alternatives}.
\newblock \emph{The Annals of Mathematical Statistics}, pages 325--346, 1968.

\bibitem[Hall and Heyde(2014)]{hall2014martingale}
P.~Hall and C.~C. Heyde.
\newblock \emph{{Martingale limit theory and its application}}.
\newblock Academic press, 2014.

\bibitem[Hallfors et~al.(2006)Hallfors, Cho, Sanchez, Khatapoush, Kim, and Bauer]{hallfors2006efficacy}
D.~Hallfors, H.~Cho, V.~Sanchez, S.~Khatapoush, H.~M. Kim, and D.~Bauer.
\newblock {Efficacy vs effectiveness trial results of an indicated “model” substance abuse program: implications for public health}.
\newblock \emph{American journal of public health}, 96\penalty0 (12):\penalty0 2254--2259, 2006.

\bibitem[Harshaw et~al.(2023)Harshaw, S{\"a}vje, Eisenstat, Mirrokni, and Pouget-Abadie]{harshaw2023design}
C.~Harshaw, F.~S{\"a}vje, D.~Eisenstat, V.~Mirrokni, and J.~Pouget-Abadie.
\newblock {Design and analysis of bipartite experiments under a linear exposure-response model}.
\newblock \emph{Electronic Journal of Statistics}, 17\penalty0 (1):\penalty0 464--518, 2023.

\bibitem[Hern{\'a}ndez et~al.(2006)Hern{\'a}ndez, Eijkemans, and Steyerberg]{hernandez2006randomized}
A.~V. Hern{\'a}ndez, M.~J. Eijkemans, and E.~W. Steyerberg.
\newblock Randomized controlled trials with time-to-event outcomes: how much does prespecified covariate adjustment increase power?
\newblock \emph{Annals of epidemiology}, 16\penalty0 (1):\penalty0 41--48, 2006.

\bibitem[Hinkelmann and Kempthorne(2007{\natexlab{a}})]{hinkelmann07}
K.~Hinkelmann and O.~Kempthorne.
\newblock \emph{{Design and Analysis of Experiments, Introduction to Experimental Design}}, volume~1.
\newblock New York: John Wiley \& Sons, 2007{\natexlab{a}}.

\bibitem[Hinkelmann and Kempthorne(2007{\natexlab{b}})]{hinkelmann2007design}
K.~Hinkelmann and O.~Kempthorne.
\newblock \emph{{Design and Analysis of Experiments, Volume 1: Introduction to Experimental Design}}, volume~1.
\newblock John Wiley \& Sons, 2007{\natexlab{b}}.

\bibitem[Ho and Chen(1978)]{ho1978lp}
S.-T. Ho and L.~H. Chen.
\newblock {An $ L_p $ Bound for the Remainder in a Combinatorial Central Limit Theorem}.
\newblock \emph{The Annals of Probability}, 6\penalty0 (2):\penalty0 231--249, 1978.

\bibitem[Hoeffding(1951)]{hoeffding1951combinatorial}
W.~Hoeffding.
\newblock {A combinatorial central limit theorem}.
\newblock \emph{The Annals of Mathematical Statistics}, pages 558--566, 1951.

\bibitem[Hu and Rosenberger(2006)]{hu2006theory}
F.~Hu and W.~F. Rosenberger.
\newblock \emph{The theory of response-adaptive randomization in clinical trials}.
\newblock John Wiley \& Sons, 2006.

\bibitem[Hudgens and Halloran(2008)]{hudgens2008toward}
M.~G. Hudgens and M.~E. Halloran.
\newblock Toward causal inference with interference.
\newblock \emph{Journal of the American Statistical Association}, 103\penalty0 (482):\penalty0 832--842, 2008.

\bibitem[Imai(2008)]{imai2008variance}
K.~Imai.
\newblock {Variance identification and efficiency analysis in randomized experiments under the matched-pair design}.
\newblock \emph{Statistics in medicine}, 27\penalty0 (24):\penalty0 4857--4873, 2008.

\bibitem[Imai et~al.(2023)Imai, Kim, and Wang]{imai2023matching}
K.~Imai, I.~S. Kim, and E.~H. Wang.
\newblock Matching methods for causal inference with time-series cross-sectional data.
\newblock \emph{American Journal of Political Science}, 67\penalty0 (3):\penalty0 587--605, 2023.

\bibitem[Imbens and Rubin(2015)]{imbens15}
G.~W. Imbens and D.~B. Rubin.
\newblock \emph{{Causal Inference in Statistics, Social, and Biomedical Sciences}}.
\newblock New York: Cambridge University Press, 2015.

\bibitem[Johansson and Schultzberg(2022)]{johansson2022rerandomization}
P.~Johansson and M.~Schultzberg.
\newblock Rerandomization: A complement or substitute for stratification in randomized experiments?
\newblock \emph{Journal of Statistical Planning and Inference}, 218:\penalty0 43--58, 2022.

\bibitem[Kempthorne(1952)]{kempthorne52}
O.~Kempthorne.
\newblock \emph{{The Design and Analysis of Experiments.}}
\newblock New York: Wiley, 1952.

\bibitem[Kempthorne(1955)]{kempthorne1955randomization}
O.~Kempthorne.
\newblock {The randomization theory of experimental inference}.
\newblock \emph{Journal of the American Statistical Association}, 50\penalty0 (271):\penalty0 946--967, 1955.

\bibitem[Laan and Robins(2003)]{laan2003unified}
M.~J. Laan and J.~M. Robins.
\newblock \emph{Unified methods for censored longitudinal data and causality}.
\newblock Springer, 2003.

\bibitem[Lee et~al.(2021)Lee, Morduch, Ravindran, Shonchoy, and Zaman]{Lee2021}
J.~N. Lee, J.~Morduch, S.~Ravindran, A.~Shonchoy, and H.~Zaman.
\newblock Poverty and migration in the digital age: Experimental evidence on mobile banking in bangladesh.
\newblock \emph{American Economic Journal: Applied Economics}, 13:\penalty0 38--71, 2021.

\bibitem[Lei and Ding(2020)]{lei2020regression}
L.~Lei and P.~Ding.
\newblock {Regression adjustment in completely randomized experiments with a diverging number of covariates}.
\newblock \emph{Biometrika}, 108\penalty0 (4):\penalty0 815--828, dec 2020.
\newblock \doi{10.1093/biomet/asaa103}.
\newblock URL \url{https://doi.org/10.1093%2Fbiomet%2Fasaa103}.

\bibitem[Leung(2022)]{leung2022causal}
M.~P. Leung.
\newblock Causal inference under approximate neighborhood interference.
\newblock \emph{Econometrica}, 90\penalty0 (1):\penalty0 267--293, 2022.

\bibitem[Li and Ding(2016)]{LD2016binary}
X.~Li and P.~Ding.
\newblock Exact confidence intervals for the average causal effect on a binary outcome.
\newblock \emph{Statistics in Medicine}, 35:\penalty0 957--960, 2016.

\bibitem[Li and Ding(2017)]{li2017general}
X.~Li and P.~Ding.
\newblock {General forms of finite population central limit theorems with applications to causal inference}.
\newblock \emph{Journal of the American Statistical Association}, 112\penalty0 (520):\penalty0 1759--1769, 2017.

\bibitem[Li and Ding(2020)]{li2020Brerandomization}
X.~Li and P.~Ding.
\newblock {Rerandomization and regression adjustment}.
\newblock \emph{Journal of the Royal Statistical Society Series B: Statistical Methodology}, 82\penalty0 (1):\penalty0 241--268, 2020.

\bibitem[Li and Small(2023)]{li2023randomization}
X.~Li and D.~S. Small.
\newblock Randomization-based test for censored outcomes: A new look at the logrank test.
\newblock \emph{Statistical Science}, 38\penalty0 (1):\penalty0 92--107, 2023.

\bibitem[Li et~al.(2018)Li, Ding, and Rubin]{li2018asymptotic}
X.~Li, P.~Ding, and D.~B. Rubin.
\newblock {Asymptotic theory of rerandomization in treatment--control experiments}.
\newblock \emph{Proceedings of the National Academy of Sciences}, 115\penalty0 (37):\penalty0 9157--9162, 2018.

\bibitem[Li et~al.(2019)Li, Ding, Lin, Yang, and Liu]{li2019peer}
X.~Li, P.~Ding, Q.~Lin, D.~Yang, and J.~S. Liu.
\newblock Randomization inference for peer effects.
\newblock \emph{Journal of the American Statistical Association}, 114:\penalty0 1651--1664, 2019.

\bibitem[Li et~al.(2020)Li, Ding, and Rubin]{li2020Arerandomization}
X.~Li, P.~Ding, and D.~Rubin.
\newblock {Rerandomization in $2^K$ factorial experiments}.
\newblock \emph{Annals of Statistics}, 48\penalty0 (1):\penalty0 43--63, 2020.

\bibitem[Lin(2013)]{lin2013agnostic}
W.~Lin.
\newblock {Agnostic notes on regression adjustments to experimental data: Reexamining Freedman’s critique}.
\newblock \emph{The Annals of Applied Statistics}, 7\penalty0 (1):\penalty0 295 -- 318, 2013.
\newblock \doi{10.1214/12-AOAS583}.
\newblock URL \url{https://doi.org/10.1214/12-AOAS583}.

\bibitem[Liu and Yang(2020)]{liu2020regression}
H.~Liu and Y.~Yang.
\newblock {Regression-adjusted average treatment effect estimates in stratified randomized experiments}.
\newblock \emph{Biometrika}, 107\penalty0 (4):\penalty0 935--948, 2020.

\bibitem[Liu et~al.(2022{\natexlab{a}})Liu, Ren, and Yang]{liu2022randomization}
H.~Liu, J.~Ren, and Y.~Yang.
\newblock {Randomization-based joint central limit theorem and efficient covariate adjustment in randomized block $2^K$ factorial experiments}.
\newblock \emph{Journal of the American Statistical Association}, pages 1--15, 2022{\natexlab{a}}.

\bibitem[Liu et~al.(2022{\natexlab{b}})Liu, Wang, and Xu]{liu2022practical}
L.~Liu, Y.~Wang, and Y.~Xu.
\newblock A practical guide to counterfactual estimators for causal inference with time-series cross-sectional data.
\newblock \emph{American Journal of Political Science}, 2022{\natexlab{b}}.

\bibitem[Lu(2016)]{lu2016covariate}
J.~Lu.
\newblock {Covariate adjustment in randomization-based causal inference for $2^K$ factorial designs}.
\newblock \emph{Statistics \& Probability Letters}, 119:\penalty0 11--20, 2016.

\bibitem[Lu and Tsiatis(2008)]{lu2008improving}
X.~Lu and A.~A. Tsiatis.
\newblock Improving the efficiency of the log-rank test using auxiliary covariates.
\newblock \emph{Biometrika}, 95\penalty0 (3):\penalty0 679--694, 2008.

\bibitem[Lu et~al.(2023{\natexlab{a}})Lu, Liu, Liu, and Ding]{lu2023design}
X.~Lu, T.~Liu, H.~Liu, and P.~Ding.
\newblock {Design-based theory for cluster rerandomization}.
\newblock \emph{Biometrika}, 110\penalty0 (2):\penalty0 467--483, 2023{\natexlab{a}}.

\bibitem[Lu et~al.(2023{\natexlab{b}})Lu, Yang, and Wang]{lu2023debiased}
X.~Lu, F.~Yang, and Y.~Wang.
\newblock Debiased regression adjustment in completely randomized experiments with moderately high-dimensional covariates.
\newblock \emph{arXiv preprint arXiv:2309.02073}, 2023{\natexlab{b}}.

\bibitem[Maclure et~al.(2006)Maclure, Nguyen, Carney, Dormuth, Roelants, Ho, and Schneeweiss]{maclure2006measuring}
M.~Maclure, A.~Nguyen, G.~Carney, C.~Dormuth, H.~Roelants, K.~Ho, and S.~Schneeweiss.
\newblock Measuring prescribing improvements in pragmatic trials of educational tools for general practitioners.
\newblock \emph{Basic \& clinical pharmacology \& toxicology}, 98\penalty0 (3):\penalty0 243--252, 2006.

\bibitem[Madow(1948)]{madow1948limiting}
W.~G. Madow.
\newblock {On the limiting distributions of estimates based on samples from finite universes}.
\newblock \emph{The Annals of Mathematical Statistics}, pages 535--545, 1948.

\bibitem[Middleton and Aronow(2015)]{middleton2015unbiased}
J.~A. Middleton and P.~M. Aronow.
\newblock {Unbiased estimation of the average treatment effect in cluster-randomized experiments}.
\newblock \emph{Statistics, Politics and Policy}, 6\penalty0 (1-2):\penalty0 39--75, 2015.

\bibitem[Moore and van~der Laan(2009)]{moore2009increasing}
K.~L. Moore and M.~J. van~der Laan.
\newblock Increasing power in randomized trials with right censored outcomes through covariate adjustment.
\newblock \emph{Journal of biopharmaceutical statistics}, 19\penalty0 (6):\penalty0 1099--1131, 2009.

\bibitem[Morgan and Rubin(2012)]{morgan2012rerandomization}
K.~L. Morgan and D.~B. Rubin.
\newblock {Rerandomization to improve covariate balance in experiments}.
\newblock \emph{Annals of Statistics}, 40\penalty0 (2):\penalty0 1263--1282, 2012.

\bibitem[Motoo(1956)]{motoo1956hoeffding}
M.~Motoo.
\newblock {On the Hoeffding’s combinatrial central limit theorem}.
\newblock \emph{Annals of the Institute of Statistical Mathematics}, 8:\penalty0 145--154, 1956.

\bibitem[Neyman(1923/1990)]{neyman1923application}
J.~Neyman.
\newblock {On the application of probability theory to agricultural experiments. Essay on principles. Section 9.}
\newblock \emph{Statistical Science}, pages 465--472, 1923/1990.

\bibitem[Neyman(1934)]{neyman1934two}
J.~Neyman.
\newblock {On the two different aspects of the representative method: the method of stratified sampling and the method of purposive selection}.
\newblock \emph{Journal of the Royal Statistical Society Series A: Statistics in Society}, 97\penalty0 (4):\penalty0 558--606, 1934.

\bibitem[Neyman and Iwaszkiewicz(1935)]{neyman1935statistical}
J.~Neyman and K.~Iwaszkiewicz.
\newblock {Statistical problems in agricultural experimentation}.
\newblock \emph{Supplement to the Journal of the Royal Statistical Society}, 2\penalty0 (2):\penalty0 107--180, 1935.

\bibitem[Noether(1949)]{noether1949theorem}
G.~E. Noether.
\newblock {On a theorem by Wald and Wolfowitz}.
\newblock \emph{The Annals of Mathematical Statistics}, 20\penalty0 (3):\penalty0 455--458, 1949.

\bibitem[Pashley and Miratrix(2021)]{pashley2021insights}
N.~E. Pashley and L.~W. Miratrix.
\newblock Insights on variance estimation for blocked and matched pairs designs.
\newblock \emph{Journal of Educational and Behavioral Statistics}, 46\penalty0 (3):\penalty0 271--296, 2021.

\bibitem[Paul and R{\'e}nyi(1959)]{erdos1959central}
E.~Paul and A.~R{\'e}nyi.
\newblock {On the central limit theorem for samples from a finite population}.
\newblock \emph{Publications of the Mathematical Institute of the Hungarian Academy of Sciences}, 4:\penalty0 49--61, 1959.

\bibitem[Petersen(1994)]{petersen1994agricultural}
R.~G. Petersen.
\newblock \emph{{Agricultural Field Experiments: Design and Analysis}}.
\newblock CRC Press, 1994.

\bibitem[Pitman(1937)]{pitman1937significance}
E.~J. Pitman.
\newblock {Significance tests which may be applied to samples from any populations}.
\newblock \emph{Supplement to the Journal of the Royal Statistical Society}, 4\penalty0 (1):\penalty0 119--130, 1937.

\bibitem[Raic(2015)]{raic2015multivariate}
M.~Raic.
\newblock {Multivariate normal approximation: Permutation statistics, local dependence and beyond}, 2015.

\bibitem[Rai{\v{c}}(2019)]{raic2019multivariate}
M.~Rai{\v{c}}.
\newblock A multivariate berry--esseen theorem with explicit constants.
\newblock 2019.

\bibitem[Reinert and R{\"o}llin(2007)]{reinert2007multivariate}
G.~Reinert and A.~R{\"o}llin.
\newblock {Multivariate normal approximation with Stein's method of exchangeable pairs under a general linearity condition}.
\newblock \emph{Annals of Probability}, 37\penalty0 (6), 2007.

\bibitem[Rigdon and Hudgens(2015)]{rigdon2015exact}
J.~Rigdon and M.~G. Hudgens.
\newblock Exact confidence intervals in the presence of interference.
\newblock \emph{Statistics and Probability Letters}, 105:\penalty0 130--135, 2015.

\bibitem[Robins(1988)]{robins1988confidence}
J.~M. Robins.
\newblock {Confidence intervals for causal parameters}.
\newblock \emph{Statistics in Medicine}, 7\penalty0 (7):\penalty0 773--785, 1988.

\bibitem[Rosenbaum(2001)]{Rosenbaum:2001}
P.~R. Rosenbaum.
\newblock Effects attributable to treatment: inference in experiments and observational studies within a discrete pivot.
\newblock \emph{Biometrika}, 88:\penalty0 219--231, 2001.

\bibitem[Rosenbaum(2002)]{rosenbaum2002observational}
P.~R. Rosenbaum.
\newblock \emph{{Observational Studies}}.
\newblock Springer-Verlag, 2002.

\bibitem[Rosenblum and van~der Laan(2010)]{Rv2010}
M.~Rosenblum and M.~J. van~der Laan.
\newblock Simple, efficient estimators of treatment effects in randomized trials using generalized linear models to leverage baseline variables.
\newblock \emph{The International Journal of Biostatistics}, 6, 2010.

\bibitem[Rubin and van~der Laan(2007)]{rubin2007doubly}
D.~Rubin and M.~J. van~der Laan.
\newblock A doubly robust censoring unbiased transformation.
\newblock \emph{The international journal of biostatistics}, 3\penalty0 (1), 2007.

\bibitem[Rubin(1974)]{rubin1974estimating}
D.~B. Rubin.
\newblock {Estimating causal effects of treatments in randomized and nonrandomized studies.}
\newblock \emph{Journal of educational Psychology}, 66\penalty0 (5):\penalty0 688, 1974.

\bibitem[Rubin(1980)]{rubin1980randomization}
D.~B. Rubin.
\newblock {Randomization analysis of experimental data: The Fisher randomization test comment}.
\newblock \emph{Journal of the American statistical association}, 75\penalty0 (371):\penalty0 591--593, 1980.

\bibitem[Rubin(1990)]{rubin1990comment}
D.~B. Rubin.
\newblock {Comment: Neyman (1923) and causal inference in experiments and observational studies}.
\newblock \emph{Statistical Science}, 5\penalty0 (4):\penalty0 472--480, 1990.

\bibitem[Rubin(2005)]{rubin2005causal}
D.~B. Rubin.
\newblock {Causal inference using potential outcomes: Design, modeling, decisions}.
\newblock \emph{Journal of the American Statistical Association}, 100\penalty0 (469):\penalty0 322--331, 2005.

\bibitem[Rubin(2008)]{rubin2008comment}
D.~B. Rubin.
\newblock {Comment: The design and analysis of gold standard randomized experiments}.
\newblock \emph{Journal of the American Statistical Association}, 103\penalty0 (484):\penalty0 1350--1353, 2008.

\bibitem[Scheff{\'e}(1959)]{scheffe1959analysis}
H.~Scheff{\'e}.
\newblock \emph{{The Analysis of Variance}}.
\newblock New York: John Wiley \& Sons, 1959.

\bibitem[Schochet et~al.(2022)Schochet, Pashley, Miratrix, and Kautz]{schochet2022design}
P.~Z. Schochet, N.~E. Pashley, L.~W. Miratrix, and T.~Kautz.
\newblock {Design-based ratio estimators and central limit theorems for clustered, blocked RCTs}.
\newblock \emph{Journal of the American Statistical Association}, 117\penalty0 (540):\penalty0 2135--2146, 2022.

\bibitem[Shi and Ding(2022)]{shi2022berry}
L.~Shi and P.~Ding.
\newblock {Berry--Esseen bounds for design-based causal inference with possibly diverging treatment levels and varying group sizes}.
\newblock \emph{arXiv preprint arXiv:2209.12345}, 2022.

\bibitem[Shi et~al.(2023)Shi, Wang, and Ding]{shi2023forward}
L.~Shi, J.~Wang, and P.~Ding.
\newblock {Forward screening and post-screening inference in factorial designs}.
\newblock \emph{arXiv preprint arXiv:2301.12045}, 2023.

\bibitem[Sj{\"o}lander et~al.(2016)Sj{\"o}lander, Frisell, Kuja-Halkola, {\"O}berg, and Zetterqvist]{sjolander2016carryover}
A.~Sj{\"o}lander, T.~Frisell, R.~Kuja-Halkola, S.~{\"O}berg, and J.~Zetterqvist.
\newblock Carryover effects in sibling comparison designs.
\newblock \emph{Epidemiology}, 27\penalty0 (6):\penalty0 852--858, 2016.

\bibitem[Splawa-Neyman(1925)]{splawa1925contributions}
J.~Splawa-Neyman.
\newblock {Contributions to the theory of small samples drawn from a finite population}.
\newblock \emph{Biometrika}, pages 472--479, 1925.

\bibitem[Sprott and Farewell(1993)]{sprott1993randomization}
D.~Sprott and V.~Farewell.
\newblock {Randomization in experimental science}.
\newblock \emph{Statistical Papers}, 34:\penalty0 89--94, 1993.

\bibitem[Stein(1972)]{stein1972bound}
C.~Stein.
\newblock {A bound for the error in the normal approximation to the distribution of a sum of dependent random variables}.
\newblock In \emph{Proceedings of the Sixth Berkeley Symposium on Mathematical Statistics and Probability, Volume 2: Probability Theory}, volume~6, pages 583--603. University of California Press, 1972.

\bibitem[Student(1923)]{student1923testing}
Student.
\newblock {On testing varieties of cereals}.
\newblock \emph{Biometrika}, pages 271--293, 1923.

\bibitem[Su and Ding(2021)]{su2021model}
F.~Su and P.~Ding.
\newblock {Model-assisted analyses of cluster-randomized experiments}.
\newblock \emph{Journal of the Royal Statistical Society Series B: Statistical Methodology}, 83\penalty0 (5):\penalty0 994--1015, 2021.

\bibitem[Su and Li(2023)]{SL2023}
Y.~Su and X.~Li.
\newblock {Treatment effect quantiles in stratified randomized experiments and matched observational studies}.
\newblock \emph{Biometrika}, 111:\penalty0 235--254, 2023.

\bibitem[Tchetgen and VanderWeele(2012)]{tchetgen2012causal}
E.~J.~T. Tchetgen and T.~J. VanderWeele.
\newblock On causal inference in the presence of interference.
\newblock \emph{Statistical methods in medical research}, 21\penalty0 (1):\penalty0 55--75, 2012.

\bibitem[Van~der Laan et~al.(2011)Van~der Laan, Rose, et~al.]{van2011targeted}
M.~J. Van~der Laan, S.~Rose, et~al.
\newblock \emph{Targeted learning: causal inference for observational and experimental data}, volume~4.
\newblock Springer, 2011.

\bibitem[VanderWeele(2015)]{vanderweele2015explanation}
T.~VanderWeele.
\newblock \emph{Explanation in causal inference: methods for mediation and interaction}.
\newblock Oxford University Press, 2015.

\bibitem[VanderWeele(2016)]{vanderweele2016mediation}
T.~J. VanderWeele.
\newblock Mediation analysis: a practitioner's guide.
\newblock \emph{Annual review of public health}, 37:\penalty0 17--32, 2016.

\bibitem[von Bahr(1976)]{von1976remainder}
B.~von Bahr.
\newblock {Remainder term estimate in a combinatorial limit theorem}.
\newblock \emph{Zeitschrift f{\"u}r Wahrscheinlichkeitstheorie und Verwandte Gebiete}, 35\penalty0 (2):\penalty0 131--139, 1976.

\bibitem[Wald and Wolfowitz(1944)]{wald1944statistical}
A.~Wald and J.~Wolfowitz.
\newblock {Statistical tests based on permutations of the observations}.
\newblock \emph{The Annals of Mathematical Statistics}, 15\penalty0 (4):\penalty0 358--372, 1944.

\bibitem[Wang et~al.(2023)Wang, Wang, and Liu]{wang2023rerandomization}
X.~Wang, T.~Wang, and H.~Liu.
\newblock {Rerandomization in stratified randomized experiments}.
\newblock \emph{Journal of the American Statistical Association}, 118\penalty0 (542):\penalty0 1295--1304, 2023.

\bibitem[Wang and Li(2022)]{wang2022rerandomization}
Y.~Wang and X.~Li.
\newblock {Rerandomization with diminishing covariate imbalance and diverging number of covariates}.
\newblock \emph{The Annals of Statistics}, 50\penalty0 (6):\penalty0 3439--3465, 2022.

\bibitem[Wang and Li(2023)]{wang2023asymptotic}
Y.~Wang and X.~Li.
\newblock Asymptotic theory of the best-choice rerandomization using the mahalanobis distance.
\newblock \emph{arXiv preprint arXiv:2312.02513}, 2023.

\bibitem[Welch(1937)]{welch1937z}
B.~L. Welch.
\newblock {On the z-test in randomized blocks and Latin squares}.
\newblock \emph{Biometrika}, 29\penalty0 (1/2):\penalty0 21--52, 1937.

\bibitem[Worrall(2010)]{worrall2010evidence}
J.~Worrall.
\newblock Evidence: philosophy of science meets medicine.
\newblock \emph{Journal of evaluation in clinical practice}, 16\penalty0 (2):\penalty0 356--362, 2010.

\bibitem[Wu and Hamada(2011)]{wu2011experiments}
C.~J. Wu and M.~S. Hamada.
\newblock \emph{{Experiments: Planning, Analysis, and Optimization}}.
\newblock John Wiley \& Sons, 2011.

\bibitem[Wu and Li(2023)]{wu2023sensitivity}
D.~Wu and X.~Li.
\newblock Sensitivity analysis for quantiles of hidden biases in matched observational studies.
\newblock \emph{arXiv preprint arXiv:2309.06459}, 2023.

\bibitem[Wu and Ding(2021)]{wu2021randomization}
J.~Wu and P.~Ding.
\newblock {Randomization tests for weak null hypotheses in randomized experiments}.
\newblock \emph{Journal of the American Statistical Association}, 116\penalty0 (536):\penalty0 1898--1913, 2021.

\bibitem[Yang and Tsiatis(2001)]{Tsiatis2001}
L.~Yang and A.~A. Tsiatis.
\newblock Efficiency study of estimators for a treatment effect in a pretest--posttest trial.
\newblock \emph{The American Statistician}, 55:\penalty0 314--321, 2001.

\bibitem[Yang et~al.(2021)Yang, Qu, and Li]{YQL2021}
Z.~Yang, T.~Qu, and X.~Li.
\newblock Rejective sampling, rerandomization, and regression adjustment in survey experiments.
\newblock \emph{Journal of the American Statistical Association}, page in press, 2021.

\bibitem[Zhang and Rosenberger(2005)]{zhang2005asymptotic}
Y.~Zhang and W.~F. Rosenberger.
\newblock On asymptotic normality of the randomization-based logrank test.
\newblock \emph{Nonparametric Statistics}, 17\penalty0 (7):\penalty0 833--839, 2005.

\bibitem[Zhao and Ding(2021)]{zhao2021covariate}
A.~Zhao and P.~Ding.
\newblock Covariate-adjusted fisher randomization tests for the average treatment effect.
\newblock \emph{Journal of Econometrics}, 225\penalty0 (2):\penalty0 278--294, 2021.

\bibitem[Zhao and Ding(2022)]{zhao2022adjust}
A.~Zhao and P.~Ding.
\newblock {To adjust or not to adjust? estimating the average treatment effect in randomized experiments with missing covariates}.
\newblock \emph{Journal of the American Statistical Association}, pages 1--11, 2022.

\bibitem[Zhao and Ding(2023)]{zhao2023covariate}
A.~Zhao and P.~Ding.
\newblock {Covariate adjustment in multiarmed, possibly factorial experiments}.
\newblock \emph{Journal of the Royal Statistical Society Series B: Statistical Methodology}, 85\penalty0 (1):\penalty0 1--23, 2023.

\bibitem[Zhao and Ding(2024)]{zhao2024no}
A.~Zhao and P.~Ding.
\newblock No star is good news: A unified look at rerandomization based on p-values from covariate balance tests.
\newblock \emph{Journal of Econometrics}, 241\penalty0 (1):\penalty0 105724, 2024.

\bibitem[Zhao et~al.(2024)Zhao, Ding, and Li]{zhao2024covariate}
A.~Zhao, P.~Ding, and F.~Li.
\newblock Covariate adjustment in randomized experiments with missing outcomes and covariates.
\newblock \emph{Biometrika}, page asae017, 2024.

\bibitem[Zhao et~al.(1997)Zhao, Bai, Chao, and Liang]{zhao1997error}
L.~Zhao, Z.~Bai, C.-C. Chao, and W.-Q. Liang.
\newblock {Error bound in a central limit theorem of double-indexed permutation statistics}.
\newblock \emph{The Annals of Statistics}, 25\penalty0 (5):\penalty0 2210--2227, 1997.

\bibitem[Zhirkov(2022)]{zhirkov2022estimating}
K.~Zhirkov.
\newblock Estimating and using individual marginal component effects from conjoint experiments.
\newblock \emph{Political Analysis}, 30\penalty0 (2):\penalty0 236--249, 2022.

\bibitem[Zigler and Papadogeorgou(2021)]{zigler2021bipartite}
C.~M. Zigler and G.~Papadogeorgou.
\newblock Bipartite causal inference with interference.
\newblock \emph{Statistical science: a review journal of the Institute of Mathematical Statistics}, 36\penalty0 (1):\penalty0 109, 2021.

\end{thebibliography}

\end{document}